\documentclass[12pt,a4paper,preprintnumbers,showkeys]{revtex4-1} 

\usepackage{tikz}
\usepackage[fleqn]{amsmath} 
\usepackage{amssymb} 
\usepackage{mathrsfs} 
\usepackage{bm}
\usepackage{graphicx} 
\usepackage{subfigure}  
\usepackage{pict2e} 
\usepackage{ulem} 
\usepackage{cases}
\usepackage[pdfstartview=FitH,
            bookmarksnumbered=true,
            bookmarksopen=true,
            colorlinks=true, 
            pdfborder=001,   
            linkcolor=blue,  
            citecolor=blue,   
            urlcolor=blue    
            ]{hyperref}      
\usepackage[left=2.5cm, right=2.5cm, top=2.5cm, bottom=2.2cm]{geometry}
\newtheorem{definition}{Definition}[section]
\newtheorem{lemma}{Lemma}[section]
\newtheorem{theorem}{Theorem}[section]
\newtheorem{corollary}{Corollary}[section]

\begin{document}


\title{The maximal spectral radius of the uniform unicyclic hypergraph with  perfect matchings}

\author{Rui Sun, Wen-Huan Wang\footnote{Corresponding author. Email: whwang@shu.edu.cn}, Zhen-Yu Ni}

\affiliation{Department of Mathematics, Shanghai University, Shanghai 200444, China}

\date{\today}

\begin{abstract}
   Let   $\mathcal{U}(n,k)$ and  $\Gamma(n,k)$ be the set of the $k$-uniform linear and nonlinear unicyclic hypergraphs having perfect matchings with  $n$ vertices respectively, where  $n\geq k(k-1)$ and $k\geq 3$.
    By using some techniques of transformations and  constructing the incidence matrices for the hypergraphs considered, we get
   the hypergraphs with the maximal spectral radii  among three kinds of hypergraphs, namely   $\mathcal{U}(n,k)$ with $n= 2k(k-1)$  and $n\geq 9k(k-1)$,
    $\Gamma(n,k)$ with $n\geq k(k-1)$,  and $\mathcal{U}(n,k)\cup \Gamma(n,k)$ with $n\geq 2k(k-1)$, where $k\geq 3$.

\end{abstract}



\maketitle


\section{Introduction}
  Let $G=(V,E)$ be a simple (i.e., no loops or multiple edges) hypergraph, where   $V=V(G)=\{v_1,v_2, \cdots, v_n\}$ is the vertex set  and $E=E(G)=\{e_1,e_2,\cdots,e_a\}$ is  the edge set with $e_i\subseteq V(G)$ for $i=1,\cdots,a$.
  $e_i$ with $1\leq i\leq a$ is  called an edge of $G$. If $|e_i|=k$ for $1\leq i\leq a$, then $G$ is called a $k$-uniform hypergraph.
 If any two edges of $G$  intersect on at most one common vertex, then
 $G$ is called a linear hypergraph.
       Let $u,v\in V(G)$.
       A path between $u$ and $v$ is denoted by $P=(v_{1},e_{1},v_{2},\ldots,v_{p},e_{p},v_{p+1})$, where $v_{1}=u$, $v_{p+1}=v$,
     all $v_{i}$ and all $e_{i}$ are distinct, and $v_i,v_{i+1}\in e_i$ for $1\leq i\leq p$.
       We say that $u$ and $v$  are connected if there exists a path of $G$ between them.
             A  hypergraph $G$ is  connected if every pair of vertices in $V(G)$ is connected.
     For $p\geq 3$, a cycle of length $p$ of  $G$ is obtained from a path $P$ of
     length $p$  by identifying $v_{1}$ with  $v_{p+1}$.

  For a $k$-uniform hypergraph $G$, if $a(k-1)-n+\omega(G)=r(G)$,
 then we call $G$ an $r(G)$-cyclic hypergraph,
 where $a$ is the number of edges in $G$,
  $n$  the number of vertices of $G$,
  $\omega(G)$  the number of components of $G$,
  and $r(G)$   the number of cycles of $G$.
  If $r(G)=1$, then  $G$ is a unicyclic hypergraph.
   In this paper, we  consider $k$-uniform connected linear and nonlinear unicyclic hypergraphs.
  %


 For  $u,v\in V(G)$ and $e\in E(G)$,  if  $\{u,v\}\subseteq e$,
  then we say that  $u$ and $v$ are adjacent and $v$ is incident with $e$.
 We denote by   $d_{G}(v)$ the degree  of $v$. Namely  $d_{G}(v)$  is the number of the edges in $G$ incident with $v$. 
 If $d_{G}(v)=1$, then we call $v$  a core vertex. If   $d_{G}(v)\geq2$, then we say that $v$ is an intersection vertex.
  For $e=\{v_1,\ldots,v_r\}\in E(G)$, if $d_{G}(v_1)\geq2$ and  $d_{G}(v_i)=1$ for $2\leq i\leq r$, then   $e$ is called a  pendent edge at $v_1$ of $G$.





 Let $\mathbb{R}$ and  $\mathbb{C}$ be the sets of real and complex numbers, respectively.
 A real tensor (or hypermatrix) $\bm{\mathcal{A}}=(a_{i_1i_2\cdots i_r})$ of $r$-order and $n$-dimension  is a multi-dimensional array
  with entries $a_{i_1i_2\cdots i_r}$ such that $a_{i_1i_2\cdots i_r}\in \mathbb{R}$, where $i_1,i_2,\cdots, i_r\in [n]$ with $[n]=\{1,2,\cdots, n\}$.
  In 2005,  Qi \cite{2005-Qi-p1302} and Lim \cite{2005-Lim-p129} independently introduced
 the concept of tensor eigenvalues and the spectra of tensors as follows.  Let $ \bm{x}=(x_1,x_2,\ldots,x_n)^{\textrm{T}}\in \mathbb{C}^n$ be an  $n$-dimensional complex column vector. Let $ \bm{x}^{[r]}=(x^{r-1}_1,x^{r-1}_2,\cdots,x^{r-1}_n)^{\textrm{T}}$, where $r$ is a positive integer. Then
  $\bm{\mathcal{A}}\bm{x}$ is a vector in $\mathbb{C}^n$ whose $i$-th component is given by
  \begin{align}\label{1}
  (\bm{\mathcal{A}}\bm{x})_i=\sum^n_{i_2,\ldots,i_r=1}a_{ii_2\cdots i_r}x_{i_2}\cdots x_{i_r},~\mbox{for~each}~i\in [n].
  \end{align}
       If there exists a number $\lambda \in \mathbb{C}$ and a nonzero vector $\bm{x} \in \mathbb{C}^n$ such that
  $\bm{\mathcal{A}}\bm{x}=\lambda \bm{x}^{[r-1]}$,
 then $\lambda$ is called an eigenvalue of $\bm{\mathcal{A}}$ and $\bm{x}$ is called an eigenvector of $\bm{\mathcal{A}}$
 corresponding to the eigenvalue $\lambda$.
 The spectral radius of $\bm{\mathcal{A}}$ is the largest modulus of the eigenvalues of $\bm{\mathcal{A}}$, i.e.,
 $\rho(\bm{\mathcal{A}})=\mbox{max}\{|\lambda| \big| ~\lambda~ \mbox{is  an  eigenvalue of } \bm{\mathcal{A}}\}$.

 For a hypergraph  $G$, there are a few tensors associated with $G$.
  The most important tensor associated with $G$ may be the adjacency tensor which is proposed by   Cooper and Dutle \cite{Cooper2012h}  in 2012 as follows. Let $G$ be a $k$-uniform hypergraph with $n$ vertices.
   The adjacency tensor of $G$ is the $k$-order  and $n$-dimension  adjacency tensor $\bm{\mathcal{A}}(G)=(a_{i_1i_2\cdots i_k})$ whose $(i_1i_2\cdots i_k)$-entry is
  \begin{equation}\label{a}
   a_{i_1i_2\cdots i_k}= \left\{
  \begin{array}{ll}
  \dfrac{1}{(k-1)!},~~~~&\mbox{if}~ \{i_1,i_2,\cdots, i_k\}\in E(G), \\
  0,                   & \mbox{otherwise}.
  \end{array}
  \right.
  \end{equation}
   The spectral radius of $\bm{\mathcal{A}}(G)$ of a $k$-uniform hypergraph $G$ is called the spectral radius of $G$ and denoted by  $\rho (G)$.
    For a  $k$-uniform hypergraph $G$ with $n$ vertices, if $G$   is connected,  then there exists a unique positive eigenvector $\bm{x}=(x_{1},\ldots,x_{n})^{\textrm{T}}$ corresponding to $\rho(G)$ with  $\Sigma_{i=1}^{n}x^{k}_{i}=1$ \cite{2013-Friedland-p738,-N-p}.
    Such positive eigenvector is called the principal eigenvector of $G$ \cite{2013-Friedland-p738,-N-p}.
    In this paper, we will consider principal eigenvector $\bm{x}$ as a mapping $x$: $V(G)\rightarrow \mathbb{R}^n$.
    The principal eigenvector $\bm{x}$ plays a key role in the spectral hypergraph theory.





   The research on the spectra of hypergraphs via tensors has drawn increasingly extensive interest.
   In recent years,  many interesting results  about the characterization of the $k$-uniform hypergraphs with the extremal spectral radii  have been obtained.
    Xiao and Wang \cite{2019-Xiao-p1392}  determined the unique hypergraphs with
 the maximal spectral radius among all the uniform supertrees and
 all the connected uniform unicyclic hypergraphs with given number  of pendent edges.
   Fan et al. \cite{2016-YiZhengFan-p845}  characterized the hypergraph(s) with the maximal spectral radius over all unicyclic hypergraphs, linear or power unicyclic hypergraphs with a given  girth, and linear or power bicyclic hypergraphs.
    Kang et al. \cite{2018-LiYingKang-p661}
    obtained  the hypergraph with the maximal spectral radius among the linear bicyclic uniform hypergraphs.
  Ouyang   et al. \cite{2017-ChenOuyang-p141} deduced the first five hypergraphs with  the maximal spectral radii among all unicyclic hypergraphs and the first three ones over all bicyclic hypergraphs.
  Among the set of supertrees \cite{2016-Li-p741,2016-Yuan-p206,2020-Wen-HuanWang-FMC} and the set of supertrees with given parameters,
  such as   a fixed diameter \cite{Xiao2018i},
  a given degree sequence \cite{2017-PengXiao-p33},
  a perfect matching \cite{2018-Zhang-p1489},
  a given number of pendent vertices  \cite{2019-Zhang-p1062},
   a given size of matching \cite{--p,2018-HaiYanGuo-p236},
   and two vertices of maximum degree \cite{2020-Wang-p}, etc,  the hypergraphs with the extremal spectral radii were also characterized.

  Motivating by the preceding works on the hypergraphs with the extremal spectral radii,
  in this paper, we consider the hypergraph with the maximal spectral radius among the set of  the $k$-uniform unicyclic hypergraphs having perfect matchings.

  A $k$-matching  of $G$ is a union of   $k$ independent edges in  $G$, where $k\geq 0$.
  A perfect  matching of $G$ is a matching that covers $V(G)$. Namely, a set $\{S_1,S_2\cdots,S_h\}$ of pairwise vertex disjoint edges of $G$ with $V(G)=S_1\cup S_2\cup \cdots\cup S_h$ is called a perfect matching of $G$.
  The hypergraphs with perfect matchings are  important in graph theory and were widely studied for their various properties \cite{2012-Treglown-p1500,2013-Khan-p1021,2016-Khan-p333}.


 Let  $G$ be a $k$-uniform unicyclic hypergraph having perfect matchings.
 Let  $Q(G)=E(G)-M(G)$, where $M(G)$ is
 the perfect matching of $G$.
  Let $\widehat{G}$ be the hypergraph induced by
 $Q(G)$, that is, $\widehat{G} = G - M(G) - S_0$, where $S_0$ is
 the set of singletons in $G - M(G)$.
 We call $\widehat{G}$ the
 capped hypergraph of $G$ and $G$ the original hypergraph of $\widehat{G}$.
 Let $|M(G)|$ and $|Q(G)|$ be the numbers of edges in $M(G)$ and
 $Q(G)$ respectively.

 Let  $\mathcal{U}(n,k)$ be the set of the $k$-uniform linear unicyclic hypergraphs having perfect matchings with
 $n$ vertices, where $k\geq 3$. 
 Let $G$ be an arbitrary hypergraph in $\mathcal{U}(n,k)$.
 Since each vertex of $G$ is saturated, we have $|M(G)|=\frac{n}{k}$, where $n$ is divisible by $k$ and $k\geq 3$.
 Thus, it follows from $n=|E(G)|(k-1)$ that
 $|Q(G)|=|E(G)|-\frac{n}{k}=\frac{n}{k(k-1)}$, where $n$ is divisible by $(k-1)k$.
 For simplicity,  we let $|Q(G)|=m$. Namely, $m$  is the number of the edges of $\widehat{G}$.
  Thus, in $\mathcal{U}(n,k)$, we get  $n=mk(k-1)$,  where  $m$ is an integer not less than 2 and  $k\geq 3$.

   Let  $\Gamma(n,k)$ be the set of the $k$-uniform nonlinear unicyclic hypergraphs having perfect matchings with  $n$ vertices, where $k\geq 3$.
   Obviously,  for each $ G\in \Gamma(n,k)$, we have $n=mk(k-1)$, where $m\geq 1$ and $m$ is the number of the edges of $\widehat{G}$.

     This paper is organized as follows.
    In Section 2, relevant notations and some lemmas which are useful for subsequent  proofs are introduced.
 In  Section 3,  by using some transformations   and  constructing the incidence matrices for the hypergraphs considered,
 the hypergraph with the maximal spectral radius is derived  among  $\mathcal{U}(n,k)$ for $n\geq 2k(k-1)$ and $k\geq 3$.
 Furthermore, among  $\Gamma(n,k)$ and $\mathcal{U}(n,k)\cup\Gamma(n,k)$,  the hypergraphs with the maximal spectral radii are characterized in Sections 4 and 5 respectively,
 where $n\geq 2k(k-1)$ and $k\geq 3$.

 \section{\label{pre}Preliminaries}

 In Section 2, we introduce some relevant definitions and necessary lemmas which are useful for us to obtain the results.

  \begin{definition}\label{definition1}\cite{2016-Li-p741}
  Let $G=(V,E)$ be a hypergraph with $u\in V$ and $e_{1},\ldots,e_{r}\in E$ such that $u\not\in e_{i}$ for $i=1,\ldots,r$, where $r\geq1$.
  Suppose that $v_{i}\in e_{i}$ and write $e_{i}'=(e_{i}\setminus\{v_{i}\})\cup\{u\}(i=1,\ldots,r)$.
  The vertices $v_{1},\ldots,v_{r}$ need not be distinct.
   Let $G'=(V,E')$ be the hypergraph with $E'=(E\setminus\{e_{1},\ldots,e_{r}\})\cup\{e_{1}',\ldots,e_{r}'\}$.
   Then we say that $G'$ is obtained from $G$ by moving edges $(e_{1},\ldots,e_{r})$ from $(v_{1},\ldots,v_{r})$ to $u$.
\end{definition}

 In $G$, if there exist two edges (denoted by $e$ and $e'$) such that $e$ and $e'$ have the same vertices,
 then we say that $e$ and $e'$ are  two multiple edges.

 \begin{lemma}\label{lemma2.1}\cite{2016-Li-p741}
  Let  $G$ and $G'$ be the two  connected hypergraphs as defined in  Definition \ref{definition1}. 
   Suppose that $G'$ contains no multiple edges. If $\bm{x}$ is the principal eigenvector of $\mathcal{A}(G)$ corresponding to $\rho(\mathcal{A}(G))$ and  $x_{u}\geq\max_{1\leq i\leq r}\{x_{v_{i}}\}$, then $\rho(\mathcal{A}(G'))>\rho(\mathcal{A}(G))$.
 \end{lemma}

 \begin{lemma}\label{lemma2.2}\cite{2017-ChenOuyang-p141}
 Let $G$ be a connected $k$-uniform hypergraph having two adjacent vertices $u_{1}$ and $u_{2}$.
 Let $G'$ be the hypergraph obtained from $G$ by moving all incident edges of $u_{2}$ (except for all common edges shared by $u_{1}$ and $u_{2}$)    from $u_{2}$ to $u_{1}$. If $G'\ncong G$, then $\rho(G')>\rho(G)$.
 \end{lemma}

  Li et al. \cite{2016-Li-p741} proposed the edge-releasing operation for $k$-uniform linear hypergraph.
   In this paper, we generalize the edge-releasing operation to $k$-uniform  hypergraph, which is shown in Definition \ref{definition2}.

 \begin{definition}\label{definition2}
 Let $G$ be a $k$-uniform hypergraph. \textcolor{black}{Let $e$ be a non-pendent edge of $G$
and $\{e_{1},\cdots,e_{r}\}$ be all the edges of $G$ adjacent to $e$, where $e_i$ and $e$ share a common vertex which is denoted by $v_i$ $(i=1,\ldots,r)$.
 Let $u$ be an arbitrary vertex of $e$.} 
  Let $G'$ be the hypergraph obtained from $G$ by moving edges $(e_{1},\ldots,e_{r})$
   \textcolor{black}{(except for all the edges which are incident with $u$)}
   from $(v_{1},\cdots,v_{r})$
   \textcolor{black}{(except for $u$)}
   to $u$. Then $G'$ is said to be obtained from $G$ by edge-releasing operation on $e$ at $u$.
\end{definition}

%

 \begin{lemma}\label{lemma2.3}
  \textcolor{black}{Let  $G$ and $G'$ be the two  connected hypergraphs as defined in Definition \ref{definition2}.}
  If $G'$ does not have multiple edges, then we have $\rho(G')>\rho(G)$.
 \end{lemma}

 \noindent\textbf{Proof}: Let $G$  and  $G'$ be the two hypergraphs as defined  in
 \textcolor{black}{Definition \ref{definition2}.} 
 Since  $e$ is a non-pendent edge of $G$, there exist some vertices in  $e$  which have degrees not less than 2.
 We denote these vertices by $v_{1},\cdots,v_{r}$, where $2\leq r\leq k$.
 By repeatedly using Lemma \ref{lemma2.2}, we get $\rho(G')>\rho(G)$.
 ~~$\Box$

 %
%
 Let $B_{G}$ be a weighted incidence matrix of a hypergraph $G$.
 We denote by $B_{G}(v,e)$ the entry of $B_{G}$ corresponding to $v$ and $e$, where $v\in V(G)$ and $e\in E(G)$.

 \begin{definition}\label{definition3}\cite{2016-LinYuanLu-p206}
 A weighted incidence matrix $B_{G}$ of a hypergraph $G$ is a $\vert V\vert\times\vert E\vert$ matrix such that for any vertex $v\in V(G)$ and any edge $e\in E(G)$, the entry $B_{G}(v,e)>0$ if $v\in e$ and $B_{G}(v,e)=0$ if $v\not\in e$.
 \end{definition}

 Let $E_{G}(v)$ be the set of the edges which are incidence with $v$, where $v\in V(G)$.

 \begin{definition}\label{definition4}\cite{2016-LinYuanLu-p206}
 A hypergraph $G$ is  $\alpha$-normal if there exists a weighted incidence matrix $B_{G}$  satisfying

 (i). $\sum_{e:e\in E_{G}(v)}B_{G}(v,e)=1$, for any $v\in V(G)$.

 (ii). $\prod_{v:v\in e}B_{G}(v,e)=\alpha$, for any $e\in E(G)$.

 Moreover, the incidence matrix $B_{G}$ is called consistent if for any cycle $v_{0}e_{1}v_{1}\cdots v_{l}$ $(v_{l}=v_{0})$
 $$\prod_{i=1}^{l}\frac{B_{G}(v_{i},e_{i})}{B_{G}(v_{i-1},e_{i})}=1.$$
 In this case, we say $G$ is consistently $\alpha$-normal.
 \end{definition}

  \begin{definition}\label{definition5}\cite{2016-LinYuanLu-p206}
 A hypergraph $G$ is  $\alpha$-subnormal if there exists a weighted incidence matrix $B_{G}$ satisfying

(i). $\sum_{e:e\in E_{G}(v)}B_{G}(v,e)\leq1$, for any $v\in V(G)$.

(ii). $\prod_{v:v\in e}B_{G}(v,e)\geq\alpha$, for any $e\in E(G)$.

 Moreover, $G$ is  strictly $\alpha$-subnormal if it is $\alpha$-subnormal but not $\alpha$-normal.
\end{definition}

 \begin{lemma}\label{lemma4}\cite{2016-LinYuanLu-p206}
Let $G$ be a connected $k$-uniform hypergraph.

 (i) $G$ is consistently $\alpha$-normal if and only if $\rho(G)=\alpha^{-\frac{1}{k}}$.

 (ii) If $G$ is $\alpha$-subnormal, then $\rho(G)\leq\alpha^{-\frac{1}{k}}$.
 Moreover, if $G$ is strictly $\alpha$-subnormal, then $\rho(G)<\alpha^{-\frac{1}{k}}$.
\end{lemma}

\section{The hypergraph with the maximal spectral radius among  $\mathcal{U}(n,k)$}

 In Section 3, we will deduce the hypergraph with the maximal spectral radius among $\mathcal{U}(n,k)$, where $n=mk(k-1)$, $m\geq2$ and $k\geq 3$.
 Some definitions are given first.

 Let $\mathcal{U}(n,k,l)$ be a subset of $\mathcal{U}(n,k)$ in which each hypergraph has a cycle  $C_l$, where $l$ is an integer with $l\geq 3$.
 Let $C_{l}=v_{1}e_{1}v_{2}e_{2}v_{3}\cdots v_{l}e_{l}v_{1}$,
where  $e_{i}=\{v_{i},v_{i,1},\ldots,v_{i,k-2},v_{i+1}\}$ with  $1\leq i\leq l-1$ and $e_{l}=\{v_{l},v_{l,1},\ldots,v_{l,k-2},v_{1}\}$.
 Let $G\in \mathcal{U}(n,k,l)$ and $M(G)$ be the perfect matching of $G$.
 According to the fact whether $C_l$ of   $G$ has at least one perfect matching edge or not,
 we classify  $\mathcal{U}(n,k,l)$ into two subsets which are denoted by $\mathcal{U}_{1}(n,k,l)$ and $\mathcal{U}_{2}(n,k,l)$,
 where $ \mathcal{U}_{1}(n,k,l)$ satisfies that each hypergraph $G$ in  it has no perfect matching edges on $C_l$ of $G$
 and $ \mathcal{U}_{2}(n,k,l)$ satisfies that each hypergraph $G$ in it has at least one perfect matching edge on $C_l$ of $G$.
 Obviously, $\mathcal{U}(n,k)=\bigcup_{l\geq3}\big(\mathcal{U}_{1}(n,k,l)\cup \mathcal{U}_{2}(n,k,l)\big)$.


 Let $\bar{\mathcal{U}}_{1}(n,k,3)$  be a subset of $\mathcal{U}_{1}(n,k,3)$  in which each hypergraph satisfies  two conditions:
    (i) each vertex in  $C_3$ must be attached by a pendent edge;  and
   (ii) at most one of the vertices in $\{v_{1},v_{2},v_{3}\}$  is attached by a hypertree which has at least $k$ edges, where $k\geq 3$.

 Let $\bar{\mathcal{U}}_{2}(n,k,3)$ be a subset of $\mathcal{U}_{2}(n,k,3)$ in which each hypergraph satisfies three conditions:
 (i) each vertex in $e_{1}\setminus\{v_{1},v_{2}\}$ of $C_3$ is a core vertex;
(ii) each vertex in $(e_{2}\cup e_{3})\setminus\{v_{1},v_{2}\}$ of  $C_3$ must be attached by a pendent edge; and
 (iii) at most one of the vertices in $\{v_{1},v_{2},v_{3}\}$ of $C_3$ is attached by a hypertree which has at least $k\geq 3$  edges,
    and each vertex in $(e_{2}\cup e_{3})\setminus\{v_{1},v_{2},v_{3}\}$ is not attached by a hypertree  which has at least $k\geq 3$ edges.

 Let $\bar{\mathcal{U}}_{2,1}(n,k,3)$ (respectively $\bar{\mathcal{U}}_{2,2}(n,k,3)$) be a subset of  $\bar{\mathcal{U}}_{2}(n,k,3)$ in which each hypergraph satisfies all  the conditions for  $\bar{\mathcal{U}}_{2}(n,k,3)$ and further satisfies that  $v_{1}$ or  $v_{2}$ (respectively $v_{3}$) of $C_3$ of $G$ is attached by a hypertree which has at least $k$ edges, where $k\geq 3$.

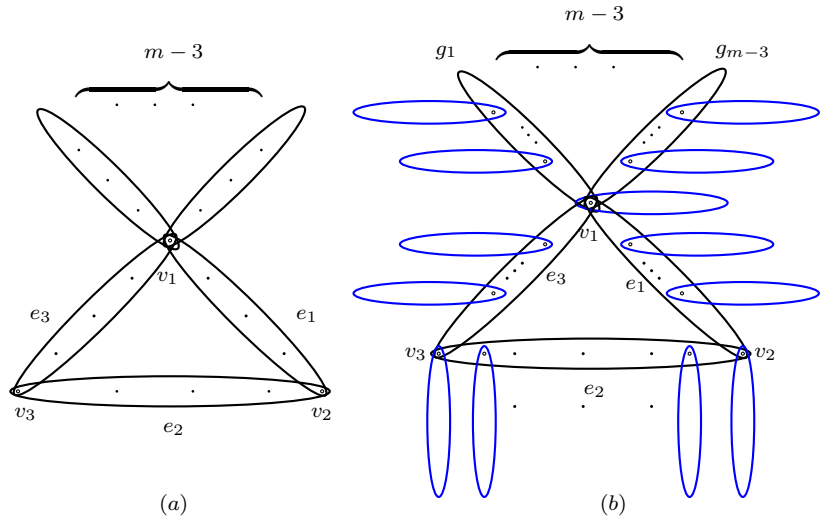
\begin{figure}
\thicklines
\begin{tikzpicture}
\draw[thick,black,rotate=45](2.4,-1.2) ellipse(1.3 and 0.2); 
\draw[thick,black,rotate=135](-0.05,-1.2) ellipse(1.3 and 0.2); 
\draw[thick,black,rotate=45](-0.19,-1.2) ellipse(1.5 and 0.2); 
\draw[thick,black,rotate=135](-2.6,-1.2) ellipse(1.5 and 0.2); 
\draw[thick,black](1.7,-2) ellipse(2.1 and 0.2); 
\draw[black](1.7,0)  circle(0.02);  
\fill(1.0,1.8) circle(0.5pt);
\fill(1.5,1.8) circle(0.5pt);
\fill(2.0, 1.8) circle(0.5pt); 
\fill(0.5, 1.2) circle(0.5pt);
\fill(0.88, 0.8) circle(0.5pt);
\fill(1.28, 0.4) circle(0.5pt);
\fill(2.1, 0.4) circle(0.5pt);
\fill(2.5, 0.8) circle(0.5pt);
\fill(2.9, 1.2) circle(0.5pt);
\draw[black](-0.3,-2) circle(0.02);  
\fill(0.2,-1.5) circle(0.02); 
\fill(0.7,-1) circle(0.02);  
\fill(1.2,-0.5) circle(0.02);  
\draw[black](3.7,-2) circle(0.02);  
\fill(3.2,-1.5) circle(0.02); 
\fill(2.7,-1) circle(0.02);  
\fill(2.2,-0.5) circle(0.02);  
\fill(2, -2) circle(0.5pt);
\fill(3,-2) circle(0.5pt);
\fill(1,-2) circle(0.5pt); 
\draw(1.75,2.5) node{{\scriptsize $m-3$}};
 \draw(1.68,2) node{$\overbrace{~~~~~~~~~~~~~~~~~~~}$};
\draw(1.75,-2.5) node{{\scriptsize $e_2$}};
\draw(3.5,-1) node{{\scriptsize $e_1$}};
\draw(0.0,-1) node{{\scriptsize $e_3$}};
\draw(1.67,-0.5) node{{\scriptsize $v_1$}};
\draw(3.7,-2.3) node{{\scriptsize $v_2$}};
\draw(-0.2,-2.3) node{{\scriptsize $v_3$}};
\draw(1.75,-3.5) node{{\scriptsize $(a)$}};
\end{tikzpicture}
\begin{tikzpicture}
\draw[thick,black,rotate=45](2.4,-1.2) ellipse(1.3 and 0.2); 
\draw[thick,black,rotate=135](-0.05,-1.2) ellipse(1.3 and 0.2); 
\draw[thick,blue](-0.42,1.2) ellipse(1.0 and 0.15); 
\draw[black](0.42,1.2) circle(0.02);  
\draw[thick,blue](3.7,1.2)  ellipse(1.0 and 0.15); 
\draw[black](2.9,1.2) circle(0.02);  
\draw[thick,blue](0.19,0.55) ellipse(1.0 and 0.15); 
\draw[black](1.1,0.55) circle(0.02);  
\draw[thick,blue](3.1,0.55)  ellipse(1.0 and 0.15); 
\draw[black](2.22,0.55) circle(0.02);  
\draw[thick,blue](2.5,0)  ellipse(1 and 0.15); 
\draw[thick,black,rotate=45](-0.19,-1.2) ellipse(1.5 and 0.2); 
\draw[thick,blue](-0.42,-1.2) ellipse(1.0 and 0.15); 
\draw[black](0.42,-1.2) circle(0.02);  
\draw[thick,blue](0.19,-0.55) ellipse(1.0 and 0.15); 
\draw[black](1.1,-0.55) circle(0.02);  
\draw[thick,black,rotate=135](-2.6,-1.2) ellipse(1.5 and 0.2); 
\draw[thick,blue](3.7,-1.2)  ellipse(1.0 and 0.15); 
\draw[black](2.9,-1.2) circle(0.02);  
\draw[thick,blue](3.1,-0.55)  ellipse(1.0 and 0.15); 
\draw[black](2.22,-0.55) circle(0.02);  
\draw[thick,black](1.7,-2) ellipse(2.1 and 0.2); 
\draw[black](1.7,0)  circle(0.02);  
\fill(1.0,1.8) circle(0.5pt);
\fill(1.5,1.8) circle(0.5pt);
\fill(2.0, 1.8) circle(0.5pt); 
\fill(0.8, 1) circle(0.5pt);
\fill(0.88, 0.9) circle(0.5pt);
\fill(0.98, 0.8) circle(0.5pt);
\fill(2.4, 0.8) circle(0.5pt);
\fill(2.5, 0.9) circle(0.5pt);
\fill(2.6, 1.0) circle(0.5pt);
\draw[black](-0.3,-2) circle(0.02);  
\fill(0.6,-1.) circle(0.02); 
\fill(0.7,-0.9) circle(0.02);  
\fill(0.8,-0.8) circle(0.02);  
\draw[thick,blue,rotate=270](2.9,-0.3) ellipse(1 and 0.15); 
\draw[thick,blue,rotate=270](2.9,0.3) ellipse(1 and 0.15); 
\draw[black](0.3,-2) circle(0.02);  
\fill(2.5,-2) circle(0.02);   
\fill(1.6,-2) circle(0.02);   
\fill(0.7,-2) circle(0.02);   
\draw[thick,blue,rotate=270](2.9,3) ellipse(1 and 0.15); 
\draw[black](3,-2) circle(0.02);  
\draw[thick,blue,rotate=270](2.9,3.7) ellipse(1 and 0.15); 
\fill(2.5,-2.7) circle(0.02);   
\fill(1.6,-2.7) circle(0.02);   
\fill(0.7,-2.7) circle(0.02);   
\draw[black](3.7,-2) circle(0.02);  
\fill(2.6,-1.) circle(0.02); 
\fill(2.5,-0.9) circle(0.02);  
\fill(2.4,-0.8) circle(0.02);  
\draw(1.75,2.5) node{{\scriptsize $m-3$}};
 \draw(1.68,2) node{$\overbrace{~~~~~~~~~~~~~~~~~~~}$};
\draw(-0.2,2) node{{\scriptsize $g_1$}};
\draw(3.7,2) node{{\scriptsize $g_{m-3}$}};
\draw(1.75,-2.5) node{{\scriptsize $e_2$}};
\draw(2.3,-1.1) node{{\scriptsize $e_1$}};
\draw(1.25,-1) node{{\scriptsize $e_3$}};
\draw(1.68,-0.43) node{{\scriptsize $v_1$}};
\draw(4,-2) node{{\scriptsize $v_2$}};
\draw(-0.6,-2) node{{\scriptsize $v_3$}};
\draw(2,-4) node{{\scriptsize $(b)$}};
\end{tikzpicture}
\caption{\label{fig-dmr} (a) $ C_{3}(v_{1},u_{0})S_{m-3,k}$ and  (b) $A_{n,k}$}
\end{figure}

  We denote by  $S_{a,k}$ the $k$-uniform linear supertree obtained from a vertex $u_0$ by attaching $a$ edges with $k$ vertices at $u_0$, where  $a\geq 1$.
   Namely, in $S_{a,k}$,  all the $a$ edges share a common vertex $u_0$.
   Let $G$ and $H$ be two hypergraphs whose vertex sets are disjoint with $v\in V(G)$ and $w\in V(H)$.
 We use   $G(v,w)H$ to denote the hypergraph obtained by identifying the vertices $v$ and $w$.
 For example,  $C_{3}(v_{1},u_{0})S_{m-3,k}$ is shown in Fig.\ \ref{fig-dmr}(a), where $C_{3}=v_1e_1v_2e_2v_3e_3v_1$ is a cycle of length 3.


  Let $A_{n,k}$ be the $k$-uniform linear unicyclic hypergraph  obtained from $C_{3}(v_{1},u_{0})S_{m-3,k}$ by attaching one pendent edge at each vertex of $C_{3}(v_{1},u_{0})S_{m-3,k}$, where $n=mk(k-1)$ and $m,k\geq 3$.
  $A_{n,k}$ is shown in Fig.\ \ref{fig-dmr}(b).

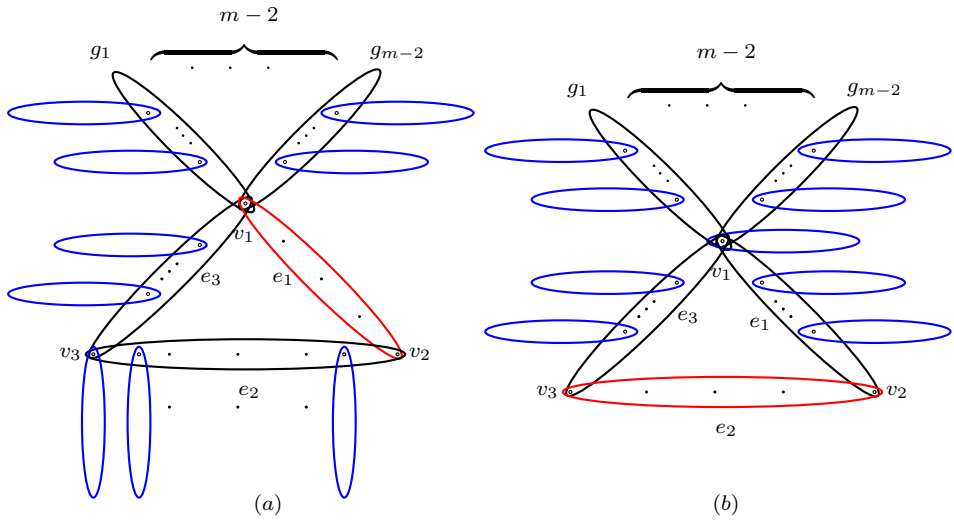
\begin{figure}
\thicklines
\begin{tikzpicture}
\draw[thick,black,rotate=45](2.4,-1.2) ellipse(1.3 and 0.2); 
\draw[thick,black,rotate=135](-0.05,-1.2) ellipse(1.3 and 0.2);  
\draw[thick,blue](-0.42,1.2) ellipse(1.0 and 0.15); 
\draw[black](0.42,1.2) circle(0.02);  
\draw[thick,blue](3.7,1.2)  ellipse(1.0 and 0.15); 
\draw[black](2.9,1.2) circle(0.02);  
\draw[thick,blue](0.19,0.55) ellipse(1.0 and 0.15); 
\draw[black](1.1,0.55) circle(0.02);  
\draw[thick,blue](3.1,0.55)  ellipse(1.0 and 0.15); 
\draw[black](2.22,0.55) circle(0.02);  
\draw[thick,black,rotate=45](-0.19,-1.2) ellipse(1.5 and 0.2); 
\draw[thick,blue](-0.42,-1.2) ellipse(1.0 and 0.15); 
\draw[black](0.42,-1.2) circle(0.02);  
\draw[thick,blue](0.19,-0.55) ellipse(1.0 and 0.15); 
\draw[black](1.1,-0.55) circle(0.02);  
\draw[thick,red,rotate=135](-2.6,-1.2) ellipse(1.5 and 0.2); 
\draw[thick,black](1.7,-2) ellipse(2.1 and 0.2); 
\draw[black](1.7,0)  circle(0.02);  
\fill(1.0,1.8) circle(0.5pt);
\fill(1.5,1.8) circle(0.5pt);
\fill(2.0, 1.8) circle(0.5pt); 
\fill(0.8, 1) circle(0.5pt);
\fill(0.88, 0.9) circle(0.5pt);
\fill(0.98, 0.8) circle(0.5pt);
\fill(2.4, 0.8) circle(0.5pt);
\fill(2.5, 0.9) circle(0.5pt);
\fill(2.6, 1.0) circle(0.5pt);
\draw[black](-0.3,-2) circle(0.02);  
\fill(0.6,-1.) circle(0.02); 
\fill(0.7,-0.9) circle(0.02);  
\fill(0.8,-0.8) circle(0.02);  
\draw[thick,blue,rotate=270](2.9,-0.3) ellipse(1 and 0.15); 
\draw[thick,blue,rotate=270](2.9,0.3) ellipse(1 and 0.15); 
\draw[black](0.3,-2) circle(0.02);  
\fill(2.5,-2) circle(0.02);   
\fill(1.6,-2) circle(0.02);   
\fill(0.7,-2) circle(0.02);   
\draw[thick,blue,rotate=270](2.9,3) ellipse(1 and 0.15); 
\draw[black](3,-2) circle(0.02);  
\fill(2.5,-2.7) circle(0.02);   
\fill(1.6,-2.7) circle(0.02);   
\fill(0.7,-2.7) circle(0.02);   
\draw[black](3.7,-2) circle(0.02);   
\fill(3.2,-1.5) circle(0.02); 
\fill(2.7,-1) circle(0.02);  
\fill(2.2,-0.5) circle(0.02);
\draw(1.75,2.5) node{{\scriptsize $m-2$}};
 \draw(1.68,2) node{$\overbrace{~~~~~~~~~~~~~~~~~~~}$};
\draw(-0.2,2) node{{\scriptsize $g_1$}};
\draw(3.7,2) node{{\scriptsize $g_{m-2}$}};
\draw(2.2,-1) node{{\scriptsize $e_1$}};
\draw(1.75,-2.5) node{{\scriptsize $e_2$}};
\draw(1.25,-1) node{{\scriptsize $e_3$}};
\draw(1.68,-0.43) node{{\scriptsize $v_1$}};
\draw(4,-2) node{{\scriptsize $v_2$}};
\draw(-0.6,-2) node{{\scriptsize $v_3$}};
\draw(2,-4) node{{\scriptsize $(a)$}};
\end{tikzpicture}
\begin{tikzpicture}
\draw[thick,black,rotate=45](2.4,-1.2) ellipse(1.3 and 0.2); 
\draw[thick,black,rotate=135](-0.05,-1.2) ellipse(1.3 and 0.2); 
\draw[thick,blue](-0.42,1.2) ellipse(1.0 and 0.15); 
\draw[black](0.42,1.2) circle(0.02);  
\draw[thick,blue](3.7,1.2)  ellipse(1.0 and 0.15); 
\draw[black](2.9,1.2) circle(0.02);  
\draw[thick,blue](0.19,0.55) ellipse(1.0 and 0.15); 
\draw[black](1.1,0.55) circle(0.02);  
\draw[thick,blue](3.1,0.55)  ellipse(1.0 and 0.15); 
\draw[black](2.22,0.55) circle(0.02);  
\draw[thick,blue](2.5,0)  ellipse(1 and 0.15); 
\draw[thick,black,rotate=45](-0.19,-1.2) ellipse(1.5 and 0.2); 
\draw[thick,blue](-0.42,-1.2) ellipse(1.0 and 0.15); 
\draw[black](0.42,-1.2) circle(0.02);  
\draw[thick,blue](0.19,-0.55) ellipse(1.0 and 0.15); 
\draw[black](1.1,-0.55) circle(0.02);  
\draw[thick,black,rotate=135](-2.6,-1.2) ellipse(1.5 and 0.2); 
\draw[thick,blue](3.7,-1.2)  ellipse(1.0 and 0.15); 
\draw[black](2.9,-1.2) circle(0.02);  
\draw[thick,blue](3.1,-0.55)  ellipse(1.0 and 0.15); 
\draw[black](2.22,-0.55) circle(0.02);  
\draw[thick,red](1.7,-2) ellipse(2.1 and 0.2); 
\draw[black](1.7,0)  circle(0.02);  
\fill(1.0,1.8) circle(0.5pt);
\fill(1.5,1.8) circle(0.5pt);
\fill(2.0, 1.8) circle(0.5pt); %
\fill(0.8, 1) circle(0.5pt);
\fill(0.88, 0.9) circle(0.5pt);
\fill(0.98, 0.8) circle(0.5pt);
\fill(2.4, 0.8) circle(0.5pt);
\fill(2.5, 0.9) circle(0.5pt);
\fill(2.6, 1.0) circle(0.5pt);
\draw[black](-0.3,-2) circle(0.02);  
\fill(2.5,-2) circle(0.02);   
\fill(1.6,-2) circle(0.02);   
\fill(0.7,-2) circle(0.02);   
\draw[black](3.7,-2) circle(0.02);  
\fill(0.6,-1.) circle(0.02); 
\fill(0.7,-0.9) circle(0.02);  
\fill(0.8,-0.8) circle(0.02);  
\fill(2.6,-1.) circle(0.02); 
\fill(2.5,-0.9) circle(0.02);  
\fill(2.4,-0.8) circle(0.02);  
\draw(1.75,2.5) node{{\scriptsize $m-2$}};
 \draw(1.68,2) node{$\overbrace{~~~~~~~~~~~~~~~~~~~}$};
\draw(-0.2,2) node{{\scriptsize $g_1$}};
\draw(3.7,2) node{{\scriptsize $g_{m-2}$}};
\draw(1.75,-2.5) node{{\scriptsize ${e}_2$}};
\draw(2.2,-1.1) node{{\scriptsize ${e}_1$}};
\draw(1.25,-1) node{{\scriptsize ${e}_3$}};
\draw(1.69,-0.43) node{{\scriptsize $v_1$}};
\draw(4,-2) node{{\scriptsize $v_2$}};
\draw(-0.6,-2) node{{\scriptsize $v_3$}};
\draw(1.75,-3.5) node{{\scriptsize $(b)$}};
\end{tikzpicture}
\caption{\label{fig-1} (a) $B_{n,k}$ and (b) $D_{n,k}$}
\end{figure}

 Let $B_{n,k}$ (respectively $D_{n,k}$) be the $k$-uniform linear unicyclic hypergraph  obtained from $C_{3}(v_{1},u_{0})S_{m-2,k}$
 by attaching one pendent edge at each vertex of $C_{3}(v_{1},u_{0})S_{m-2,k}$ except for all the vertices which are incident with  $e_1$ (respectively $e_2$), where $n= mk(k-1)$, $m\geq 2$ and $k\geq 3$.
 $B_{n,k}$ and $D_{n,k}$ are shown in Fig.\ \ref{fig-1}.


 Obviously, we have $A_{n,k}\in \bar{\mathcal{U}}_{1}(n,k,3)$, $B_{n,k}\in\bar{\mathcal{U}}_{2,1}(n,k,3)$, and $D_{n,k}\in\bar{\mathcal{U}}_{2,2}(n,k,3)$.

 To obtain the hypergraph with the maximal spectral radius in $\mathcal{U}(n,k)$ (as shown in  Theorem \ref{theorem3.6}), several lemmas are introduced first.
 Lemmas \ref{lemmawang}--\ref{lemma33.1} are introduced to get the  hypergraph with the maximal spectral radius in $\mathcal{U}_{1}(n,k,l)$ (as shown in  Corollary \ref{corollary3.1}).
 Lemmas \ref{lemma3.1}--\ref{lemma33.2} are proposed to obtain the  hypergraph with the maximal spectral radius in $\mathcal{U}_{2}(n,k,l)$ (as shown in  Corollary \ref{corollary3.2}).

  \begin{lemma}\label{lemmawang}
  Let $G\in \mathcal{U}(n,k,l)$,
  where $n\geq 3k(k-1)$,  $k\geq3$ and $l\geq 4$.
  Let $e$ be a perfect matching edge of $G$ and
  $e$ is not a pendent edge.
  Let $G_0$ be the hypergraph obtained from  $G$ by applying the edge-releasing operation on $e$ at a vertex of $e$ such that  $e$ of $G_0$ is a pendent edge.
 We have $\rho(G_0)>\rho(G)$, where
 $G_0\in \mathcal{U}(n,k)$.
 \end{lemma}

  \noindent\textbf{Proof}: Let $n\geq 3k(k-1)$ and  $k\geq3$.  Let $G$ and $G_0$ be the two hypergraphs as defined in Lemma \ref{lemmawang}. 
  \textcolor{black}{Let $M(G)$ be the perfect matching of $G$ and $e\in M(G)$. Then}
  all the edges which are \textcolor{black}{adjacent to} $e$ are edges of $Q(G)$.
    After applying the  edge-releasing operation on $e$ at a vertex of $e$,  $G_0$ has
    \textcolor{black}{the perfect matching $M(G)$ and $G_0\in \mathcal{U}(n,k)$. Since $G$ is linear, $e$ of $G_0$ is a pendent edge.}
      By Lemma \ref{lemma2.3},
      we have $\rho(G_{0})>\rho(G)$.
      ~~$\Box$

 \begin{lemma}\label{lemma3.2}
  Let $G\in \mathcal{U}_{1}(n,k,l)$, where $n\geq 3k(k-1)$, $k\geq3$ and  $l\geq4$.
  There exists one hypergraph $G^{\star}\in \mathcal{U}_{1}(n,k,3)$ such that $\rho(G^{\star})>\rho(G)$.
 \end{lemma}

  \noindent\textbf{Proof}: Let $n\geq 3k(k-1)$, $k\geq3$ and  $l\geq4$.
  Let $G\in \mathcal{U}_{1}(n,k,l)$.
  We denote the  perfect matching of $G$ by $M(G)$.
  The cycle contained in $G$ is denoted by \textcolor{black}{$C_l=v_{1}e_{1}\ldots v_{l}e_{l}v_{1}$.}
   Bearing the definition of  $\mathcal{U}_{1}(n,k,l)$ in mind,
    $e_{i}\not\in M(G)$ and each vertex of $e_{i}$ is incident with  an edge in $M(G)$, where $1\leq i\leq l$.
    Let  $\bm{x}$ be the principal eigenvector of $G$ corresponding to $\rho(G)$. Without loss of generality, we suppose $x_{v_{1}}\geq x_{v_{2}}$.

 Let $G_{1}$ be the hypergraph obtained from $G$ by moving $e_{2}$ from $v_{2}$ to $v_{1}$.
 It is  noted that all the edges which are incident with $v_2$ of $e_1$ of $G$ (except for $e_{2}$) remain unchanged.
 Therefore, $M(G)$ is the perfect matching of $G_{1}$  and $G_{1}$ contains a cycle $C_{l-1}$. Namely $G_{1}\in \mathcal{U}_{1}(n,k,l-1)$.
  By Lemma \ref{lemma2.1}, we have  $\rho(G_{1})>\rho(G)$.
  By repeatedly using the same procedure, we finally get a hypergraph $G^{\star}\in \mathcal{U}_{1}(n,k,3)$ such that $\rho(G^{\star})>\rho(G)$.
~~$\Box$


 \begin{lemma}\label{lemma3.4}
  Let $G\in \mathcal{U}_{1}(n,k,3)$, where $n\geq 3k(k-1)$ and $k\geq3$.
   There exists a  hypergraph $G^{\diamond}\in\bar{\mathcal{U}}_{1}(n,k,3)$ such that $\rho(G^{\diamond})\geq\rho(G)$ with the equality iff  (if and only if) $G\cong G^{\diamond}$.
 \end{lemma}

 \noindent\textbf{Proof}:  Let $n=mk(k-1)$,  $m\geq 3$ and $k\geq 3$.  When $m= 3$, obviously Lemma \ref{lemma3.4} holds. Next, let $m\geq 4$.
  Let $G\in \mathcal{U}_{1}(n,k,3)$. Let $M(G)$ be the perfect matching of $G$.
  \textcolor{black}{By the definition of $\mathcal{U}_{1}(n,k,3)$,}
  each vertex in $C_3$ of $G$ is incident with an edge
  (denoted by $e$) of $M(G)$
  \textcolor{black}{and} $e$ is not an edge on $C_3$ of $G$.
  By applying the edge-releasing operation on $e$ at a vertex of $e$  and by Lemma \ref{lemmawang}, we get a hypergraph $G_{2}$ such that $\rho(G_{2})\geq \rho(G)$ with the equality iff $G\cong G_2$,
   where $G_{2}$ satisfies that (i) each vertex in  $C_3$ must be attached by a pendent edge;
  (ii)  the vertex in $C_3$   may or may not be attached by a hypertree which has at least $k\geq 3$ edges.

   Since $m\geq 4$, in $G_{2}$,  there exists at least one vertex of $C_3$  of $G_{2}$ which is attached by a hypertree having at least $k\geq 3$ edges.
   Let $\bm{x}$ be the principal eigenvector of $G_{2}$ corresponding to $\rho(G_{2})$.
   Among all the vertices of $C_3$  of $G_{2}$ which are attached by hypertrees having at least $k\geq 3$ edges,
   we assume the vertex $w$ at $C_3$ of $G_{2}$ has the largest component among the  vector $\bm{x}$ and we denote it by $x_w$.
   Without loss of generality, we suppose $w$ belongs to $e_1$ of $C_3$ of $G_{2}$.
   By Lemma \ref{lemma2.1}, we  get $\rho(G_{3})\geq \rho(G_{2})$
   with the equality iff $G_{2}\cong G_{3}$,
   where       $G_{3}$ satisfies that  each  vertex in  $C_3$  of $G_{3}$  is attached by a pendent edge,
   and only one vertex (namely $w$) at $C_3$ of $G_3$ is also attached by a hypertree having at least $k\geq 3$ edges.
    If $w$ is one of $v_1$ and $ v_2$,  then we get Lemma \ref{lemma3.4}.
    Next, we assume  $w\neq v_1,v_2$.
   Let $\bm{y}$ be the principal eigenvector of $G_3$ corresponding to $\rho(G_{3})$.
    Two cases are considered as follows.

\noindent\textbf{Case (i).} $y_{v_{1}}\geq y_{w}$.

  Let $G_{4}$ be  the hypergraph obtained from $G_3$ by moving all the edges which are incident with $w$ (except for the edge $e_1$ and the pendent edge attached at $w$) from $w$ to $v_{1}$.
  Obviously, $G_{4}\in\bar{\mathcal{U}}_{1}(n,k,3)$.
  By Lemma \ref{lemma2.1}, we obtain  $\rho(G_{4})>\rho(G_{3})$.

 \noindent\textbf{Case (ii).} $y_{v_{1}}< y_{w}$.

 Let $G_5$ be  the hypergraph obtained from $G_3$ by moving $e_{3}$ from $v_{1}$ to $w$.
 It is noted that the pendent edge attached at $v_1$ of $C_3$ of  $G_3$ remains unchanged.
  Obviously, $G_{5}\in\bar{\mathcal{U}}_{1}(n,k,3)$.
 By  Lemma \ref{lemma2.1}, we have  $\rho(G_{5})>\rho(G_{3})$.

  By combining the above proofs, we have Lemma \ref{lemma3.4}.~~$\Box$

 \begin{lemma}\label{lemma33.1}
 Let $G\in\bar{\mathcal{U}}_{1}(n,k,3)$, where $n\geq 3k(k-1)$ and $k\geq 3$.
  We have $\rho(A_{n,k})\geq\rho(G)$ with the equality iff $G\cong A_{n,k}$.
 \end{lemma}

  \noindent\textbf{Proof}:
  Let $G\in\bar{\mathcal{U}}_{1}(n,k,3)$ with $n= mk(k-1)$, $m\geq 3$ and $k\geq 3$.
   If $m=3,4$, then  $\bar{\mathcal{U}}_{1}(n,k,3)=\{ A_{n,k}\}$ and  we have Lemma \ref{lemma33.1}.
   Next, let $m\geq 5$.
   Since  $m\geq 5$,  each vertex in $C_3$ of $G$ is attached by a pendent edge and
   we assume $v_{1}$ of $C_3$ of $G$ is attached by a hypertree (denoted by $\mathcal{T}$) which has at least two edges
   belonging to  $Q(G)$ and $T$ maybe has  perfect matching edges which are not pendent edges.
 By applying the edge-releasing on an arbitrary perfect matching edge (denoted by $e$) of $T$ at a vertex of $e$ which is  not a pendent edge, it follows from  Lemma \ref{lemmawang} that we  get a hypergraph $G_{6}$ such that $\rho(G_{6})\geq\rho(G)$ with the equality iff $G\cong G_{6}$,
  where $G_{6}\in\bar{\mathcal{U}}_{1}(n,k,3)$ and all the  perfect matching edges of $G_{6}$ are pendent edges.

   If $G_{6}\cong A_{n,k}$, then we get Lemma \ref{lemma33.1}.
   Otherwise, we assume $G_{6}\ncong A_{n,k}$.
   Obviously, 
   $v_{1}$ of $C_3$ of $G_{6}$ is attached by a hypertree (denoted by $T'$) which has at least two edges
   belonging to $Q(G_{6})$.
  As $G_{6}\ncong A_{n,k}$,
  in $G_{6}$, there exists one edge (denoted by $g=\{u_{1},...,u_{k}\}$) which satisfies the following conditions:
  (i) $g$ belongs to $Q(G_{6})$;
  \textcolor{black}{(ii) $g$ is not incident with $v_1$;  and
   (iii) $v_1$ and $u_1$ are incident with a common edge.}
    Let  $\bm{x}$ be the principal eigenvector of $G_{6}$ corresponding to $\rho(G_{6})$.
    Two cases are considered as follows.

  \noindent\textbf{Case (i).} $x_{v_{1}}\geq x_{u_{1}}$.

   Let $G_{7}$ be the hypergraph obtained from $G_{6}$ by moving $g$ from $u_{1}$ to $v_{1}$.
   Obviously, $G_{7} \in\bar{\mathcal{U}}_{1}(n,k,3)$.
   By Lemma \ref{lemma2.1},  we obtain $\rho(G_{7})>\rho(G_{6})$.


  \noindent\textbf{Case (ii).}  $x_{v_{1}}<x_{u_{1}}$.

   Let $G_{8}$ be the hypergraph obtained from $G_{6}$ by moving all the edges which are incident with $v_1$
   (except for the pendent edge attached at $v_1$ and the common edge which is incident with $v_1$ and $u_1$) from $v_1$  to $u_1$.
   Obviously, $G_{8} \in\bar{\mathcal{U}}_{1}(n,k,3)$.
   By Lemma \ref{lemma2.1},  we  obtain $\rho(G_{8})>\rho(G_{6})$.

   By repeatedly using the same procedures as these in  Cases (i) and  (ii), we finally get $\rho(A_{n,k})>\rho(G_{6})$.

   By combining the above proofs, we obtain $\rho(A_{n,k})\geq\rho(G)$ with the equality iff $G\cong A_{n,k}$. ~~$\Box$

 By Lemmas \ref{lemma3.2}--\ref{lemma33.1}, we get Corollary \ref{corollary3.1} as follows.

  \begin{corollary}\label{corollary3.1}
 Let $G\in \mathcal{U}_{1}(n,k,l)$, where $n\geq 3k(k-1)$ and $k,l\geq3$.
 We have $\rho(A_{n,k})\geq\rho(G)$ with the equality iff $G\cong A_{n,k}$.
 \end{corollary}


 \begin{lemma}\label{lemma3.1}
 Let $G\in \mathcal{U}_{2}(n,k,l)$, where $n\geq 2k(k-1)$, $k\geq 3$ and $l\geq4$.
 There exists one hypergraph $\widetilde{G}\in \mathcal{U}_{1}(n,k,p)\cup \mathcal{U}_{2}(n,k,3)$ such that
 $\rho(\widetilde{G})>\rho(G)$, where $3\leq p\leq l-1$.
\end{lemma}

\noindent\textbf{Proof}:
 Let $G\in \mathcal{U}_{2}(n,k,l)$, where $l\geq4$.
Let $M(G)$ be the perfect matching of $G$. By the definition of $\mathcal{U}_{2}(n,k,l)$,
  there exists one edge (denoted by $e$) on the cycle $C_l$ of $G$ which belongs to $M(G)$. 
  Let $v$ be an arbitrary vertex of $e$.
   Let $G_{9}$ be the hypergraph obtained from $G$ by applying the edge-releasing operation on $e$ at $v$.
 It follows from Lemma \ref{lemmawang} that $\rho(G_{9})>\rho(G)$\textcolor{black}{,}
  $e$ becomes a pendent edge in $G_{9}$ and
   the  number of the perfect matching edges on the cycle contained in $G_{9}$ decreases by 1.
   Obviously, we have $G_{9}\in \mathcal{U}_2(n,k,l-1)$.
 By repeatedly using the same procedure, we finally get a hypergraph
 $\widetilde{G}\in \mathcal{U}_{1}(n,k,p)\cup \mathcal{U}_{2}(n,k,3)$ such that
 $\rho(\widetilde{G})>\rho(G)$, where $3\leq p\leq l-1$.
~~$\Box$

 \begin{lemma}\label{lemma3.3}
 Let $G\in \mathcal{U}_{2}(n,k,3)$, where $n\geq 2k(k-1)$ and $k\geq 3$.
 There exists a hypergraph $G^{\ast}\in\bar{\mathcal{U}}_{2}(n,k,3)$ such that $\rho(G^{\ast})\geq\rho(G)$
 with the equality iff  $G\cong G^{\ast}$.
\end{lemma}

\noindent\textbf{Proof}: Let  $n\geq 2k(k-1)$ and $k\geq 3$.
 Let $G\in \mathcal{U}_{2}(n,k,3)$ and $C_{3}=v_{1}e_{1}v_{2}e_{2}v_{3}e_{3}v_{1}$ be the cycle contained in $G$.
 Let $M(G)$ be the perfect matching of $G$.
 According to the definition of $\mathcal{U}_{2}(n,k,3)$, there exists one edge on $C_{3}$ which belongs to $M(G)$.
 Without loss of generality, we suppose that $e_{1}=\{v_{1},v_{1,1},\ldots,v_{1,k-2},v_{2}\}$ belongs to $M(G)$.
 If there exists one vertex in $e_{1}\setminus \{v_1,v_2\}$ such that its degree is greater than 1, then we suppose  this vertex is   $v_{1,1}$.
 Let $G_{10}$ be the hypergraph obtained from $G$ by moving all the edges (except for $e_{1}$) which are incident with  $v_{1,1}$ from  $v_{1,1}$ to  $v_{1}$.
 Obviously, $G_{10}$ has the perfect matching $M(G)$ and the number of the core vertices of $G_{10}$ increases by 1. By Lemma \ref{lemma2.2}, $\rho(G_{10})> \rho(G)$.
 By repeatedly using the same procedure, we finally get a hypergraph $G_{11}$ such that $\rho(G_{11})\geq \rho(G)$ with the equality iff $G\cong G_{11}$,
 where $G_{11}\in \mathcal{U}_{2}(n,k,3)$, $G_{11}$ has the perfect matching $M(G)$ and each vertex in $e_{1}\backslash\{v_{1},v_{2}\}$ of $G_{11}$ is a core vertex.

  Since $e_{2}$ and $e_{3}$ of $C_3$ of $G_{11}$ are the edges of $Q(G_{11})$, each vertex in $(e_{2}\cup e_{3})\setminus\{v_{1},v_{2}\}$ is incident with  an edge in $M(G)$.
  By  Lemma \ref{lemmawang}, we get a hypergraph $G_{12}$ such that $\rho(G_{12})\geq\rho(G_{11})$
   with the equality iff $G_{12}\cong G_{11}$, where $G_{12}$ satisfies the following three conditions: (i) each vertex in $e_{1}\setminus\{v_{1},v_{2}\}$ of $C_3$ is a core vertex;
(ii) each vertex in $(e_{2}\cup e_{3})\setminus\{v_{1},v_{2}\}$ of  $C_3$ must be attached by a pendent edge; and
 (iii) the vertex in $e_{2}\cup e_{3}$ of  $C_3$  may or may not be attached by a hypertree which has at least $k\geq 3$ edges.

 Furthermore, by the methods similar to those for the second paragraph, Cases (i) and (ii) in Lemma \ref{lemma3.4}, we obtain a hypergraph $G^{\ast}$ such that $\rho(G^{\ast})\geq\rho(G_{12})$ with the equality iff
 $G^{\ast} \cong G_{12}$, where $G^{\ast}\in\bar{\mathcal{U}}_{2}(n,k,3)$.
 Therefore, we get 
 \textcolor{black}{Lemma \ref{lemma3.3}.}
  ~~$\Box$

\begin{lemma}\label{lemma33.3}
  For $n\geq 3k(k-1)$ and $k\geq3$, we  have  $\rho(A_{n,k})>\rho(B_{n,k})$, where  $A_{n,k}$ and $B_{n,k}$ are the two hypergraphs as shown in
  Fig.\ \ref{fig-dmr}(b) and  Fig.\ \ref{fig-1}(a), respectively.
\end{lemma}

 \noindent\textbf{Proof}:  Let $n\geq 2k(k-1)$ and $k\geq 3$. We construct a weighted incidence matrix $B_{A_{n,k}}$ for $A_{n,k}$ as follows.
 $$
 B_{A_{n,k}}(v,e)=\left\{
 \begin{aligned}
 &0& &v\not\in e,\\
 &1& &v\in e~ and~ v\text{ is a core vertex},\\
 &\alpha & & v\in e,~ e~ is~ a ~pendent~ edge ~and~ d_{A_{n,k}}(v)=2,\\
 &1-\alpha & &v\in e,~ e~ is~not~ a ~pendent~ edge ~and~ d_{A_{n,k}}(v)=2,\\
 &\dfrac{\alpha }{(1-\alpha )^{k-1}}& &(v,e)=(v_1,g_{i}), where~ i=1,\cdots,m-3, ~ and~ m\geq3 ,\\
 &x_{0}& &(v,e)=(v_1,e_{1}),\\
 &y_{0}& &(v,e)=(v_2,e_{1}),\\
 &c_{0}& &(v,e)=(v_1,e_{3}),\\
 &d_{0}& &(v,e)=(v_2,e_{2}),\\
 &\dfrac{\alpha }{c_{0}(1-\alpha )^{k-2}}& &(v,e)=(v_3,e_{3}),\\
 &\dfrac{\alpha }{d_{0}(1-\alpha )^{k-2}}& &(v,e)=(v_3,e_{2}).\\
 \end{aligned}
 \right.
 $$
  where $x_{0}$, $y_{0}$, $c_{0}$,  $d_{0}$, and $\alpha$ satisfy the following five equations:
 \begin{numcases}{}
 x_{0}+\alpha +c_{0}+\dfrac{(m-3)\alpha }{(1-\alpha )^{k-1}}=1,&\label{1.1}  \\
 y_{0}+\alpha +d_{0}=1,& \label{1.2}\\
 \dfrac{\alpha }{c_{0}(1-\alpha )^{k-2}}+\dfrac{\alpha }{d_{0}(1-\alpha )^{k-2}}+\alpha =1,&\nonumber\\
 x_{0}y_{0}(1-\alpha )^{k-2}=\alpha ,& \label{1.3}\\
 \dfrac{x_{0}}{y_{0}}\cdot\dfrac{d_{0}^{2}}{c_{0}^{2}}=1.&\nonumber
 \end{numcases}
 We check that  $\sum_{e:e\in E_{A_{n,k}}(v)}B_{A_{n,k}}(v,e)=1$ for any $v\in V(A_{n,k})$,
  $\prod_{v:v\in e}B_{A_{n,k}}(v,e)=\alpha $ for any $e\in E(A_{n,k})$, and  $B_{A_{n,k}}$ is consistent.
   Thus, $A_{n,k}$ is consistently $\alpha$-normal.
   By Lemma \ref{lemma4}(i), we have $\rho({A_{n,k}})=\alpha^{-\frac{1}{k}}$.

We construct a weighted incidence matrix $B_{B_{n,k}}$ for $B_{n,k}$ as follows.
 Let $B_{B_{n,k}}(v,e)=0$ if $v\not\in e$;
   $B_{B_{n,k}}(v,e)=1$ if $v\in e$ and $v$ is a core vertex;
  $B_{B_{n,k}}(v,e)=\alpha$ if $v\in e$, $e$ is a pendent edge and $d_{B_{n,k}}(v)=2$;
 $B_{B_{n,k}}(v,e)=1-\alpha$ if $v\in e$, $v\neq v_2$, $e$ is not a pendent edge, and $d_{B_{n,k}}(v)=2$;
 $B_{B_{n,k}}(v_1,g_i)=\dfrac{\alpha}{(1-\alpha)^{k-1}}$ for $i=1,\cdots ,m-2$ and $m\geq 3$;
 $B_{B_{n,k}}(v_1,e_{1})=x_{1}$;
 $B_{B_{n,k}}(v_2,e_{1})=y_{1}$;
 $B_{B_{n,k}}(v_1,e_{3})=c_{0}$;
 $B_{B_{n,k}}(v_2,e_{2})=d_{0}$;
 $B_{B_{n,k}}(v_3,e_{3})=\dfrac{\alpha }{c_{0}(1-\alpha )^{k-2}}$; and
 $B_{B_{n,k}}(v_3,e_{2})=\dfrac{\alpha }{d_{0}(1-\alpha )^{k-2}}$,
%
 where $x_{1}$, $y_{1}$, $c_{0}$,  $d_{0}$, and $\alpha$ satisfy (\ref{1.4}) and (\ref{1.5}) as follows:
 \begin{numcases}{}
 x_{1}+c_{0}+\dfrac{(m-2)\alpha }{(1-\alpha )^{k-1}}=1, &\label{1.4}\\
 y_{1}+d_{0}=1. &\label{1.5}
 \end{numcases}
 We can verify that  $\sum_{e:e\in E_{B_{n,k}}(v)}B_{B_{n,k}}(v,e)=1$ for any $v\in V(B_{n,k})$ and
  $\prod_{v:v\in e}B_{B_{n,k}}(v,e)=\alpha$ for any $e\in E(B_{n,k})$ and $e\neq e_{1}$.
  Next, we prove $\prod_{v: v\in e_{1}}B_{B_{n,k}}(v,e_{1})>\alpha$.



  We get
  \begin{align}
  &\prod_{v:v\in e_{1}}B_{B_{n,k}}(v,e_{1})-\alpha \nonumber\\
  &=x_{1}y_{1}-\alpha \nonumber\\
  &=\big(x_{0}+\alpha -\dfrac{\alpha }{(1-\alpha )^{k-1}}\big)(y_{0}+\alpha )-\alpha \label{1}\\
  &=\alpha (y_{0}+\alpha -1)+x_{0}y_{0}+x_{0}\alpha -\dfrac{\alpha }{(1-\alpha )^{k-1}}(y_{0}+\alpha )\label{w1}\\
  &>\alpha (y_{0}+\alpha -1)+\dfrac{\alpha }{(1-\alpha )^{k-1}}(1-\alpha -y_{0})\label{2}\\
  &=\alpha d_{0}\big(\dfrac{1}{(1-\alpha )^{k-1}}-1\big)\label{3}\\
  &>0\label{4}
     \end{align}
 It is noted that (\ref{1}) follows from $x_{1}=x_{0}+\alpha -\dfrac{\alpha }{(1-\alpha )^{k-1}}$ (by (\ref{1.1}) and (\ref{1.4})) and
 $y_{1}=y_{0}+\alpha$ (by (\ref{1.2}) and (\ref{1.5})).
 By (\ref{1.2}) and $d_{0}>0$,  we get $y_{0}<1-\alpha $. Furthermore, by $0<y_{0}<1-\alpha <1$ and (\ref{1.3}), we obtain
  $x_{0}=\dfrac{\alpha }{y_{0}(1-\alpha )^{k-2}}>\dfrac{\alpha }{(1-\alpha )^{k-1}}$.
     Substituting this inequality and  (\ref{1.3}) into (\ref{w1}), we have  (\ref{2}).
  By (\ref{1.2}), we obtain (\ref{3}).
  Since $0<\alpha <1$, we have (\ref{4}).

%
  By the above proofs, we get that   $B_{n,k}$ is strictly $\alpha$-subnormal.
 Therefore, by Lemma \ref{lemma4}(ii), we have $\rho({B_{n,k}})<\alpha^{-\frac{1}{k}}$.
 Thus, by Lemma \ref{lemma4},  we obtain $\rho(A_{n,k})>\rho(B_{n,k})$, where $n\geq 3k(k-1)$ and $k\geq3$.
 ~~$\Box$


 \begin{lemma}\label{lemma33.4}
 Let $G\in\bar{\mathcal{U}}_{2,1}(n,k,3)$,
 where $n \geq 3k(k-1)$ and $k\geq 3$.
 We have  $\rho(A_{n,k})>\rho(G)$.
 \end{lemma}

 \noindent\textbf{Proof}: Let $n =mk(k-1)$ with $m\geq 3$ and $k\geq 3$.
 Let $G\in\bar{\mathcal{U}}_{2,1}(n,k,3)$.
 When $m=3$,  we get Lemma \ref{lemma33.4} since  $G\cong B_{n,k}$ and $\rho(A_{n,k})>\rho(B_{n,k})$ (by Lemma \ref{lemma33.3}).
  Next,  let  $m \geq 4$.

    Let  $M(G)$ be the perfect matching of $G$.
    Since $m \geq 4$, bearing the definition of $\bar{\mathcal{U}}_{2,1}(n,k,3)$ in mind,
      we assume $v_{1}$ of $C_3$ of $G$ is attached by a hypertree (denoted by $T^{\ast}$) which has at least two edges
   belonging to  $Q(G)$ and $T^{\ast}$ maybe has  perfect matching edges which are not pendent edges.
   By applying the edge-releasing on an arbitrary perfect matching
   edge (denoted by $e$) of $T^{\ast}$ at a vertex of $e$ which is
    not a pendent edge,
    it follows from  Lemma \ref{lemmawang} that    
   we can finally get a hypergraph $G_{13}$ such that $\rho(G_{13})\geq\rho(G)$ with the equality iff $G\cong G_{13}$,
   where $G_{13}\in \bar{\mathcal{U}}_{2,1}(n,k,3)$ and all the  perfect matching edges of $G_{13}$
   are pendent edges except for $e_1$ contained in $C_3$ of $G_{13}$. 

   If $G_{13}\cong B_{n,k}$, then by  Lemma \ref{lemma33.3}, we get $\rho(A_{n,k})>\rho(B_{n,k})\geq\rho(G)$. Namely, Lemma \ref{lemma33.4} holds.
   Otherwise, we assume $G_{13}\ncong B_{n,k}$.
  Since $m \geq 4$, it is noted that $v_{1}$ of $C_3$ of $G_{13}$ is attached by a hypertree (denoted by $T^{\star}$) which has at least two edges  belonging to $Q(G_{13})$.
   As $G_{13}\ncong B_{n,k}$, in $G_{13}$,
   there exists one edge (denoted by $\dot{g}=\{w_{1},...,w_{k}\}$) which satisfies the following conditions:
  (i) $\dot{g}\in E(T^{\star})$;
  (ii) $\dot{g}$ is not a pendent edge;
  (iii)  $v_1$ and $w_1$ are incident with a common edge;   and
   (iv) $\dot{g}$ is not incident with $v_1$.
    Let  $\bm{x}$ be the principal eigenvector of $G_{13}$ corresponding to $\rho(G_{13})$.
  Two cases are considered as follows.
%

\noindent\textbf{ Case (i).} $x_{v_{1}}<x_{w_{1}}$.

  Let $G_{14}$ be the hypergraph obtained from $G_{13}$ by moving all the edges which are incident with $v_1$
   (except for $e_1$ of $C_3$ of $G_{13}$ and the common edge which is incident with $v_1$ and $w_1$) from $v_1$ to $w_1$.
   Obviously, $G_{14} \in{{\mathcal{U}}}_{2}(n,k,4)$, 
   where each vertex of $e_1 \setminus\{v_{1},v_{2}\}$ of $C_4$ of $G_{14}$ is a core vertex and only $e_1$ is  the  perfect matching edge on $C_4$ of $G_{14}$.
   By Lemma \ref{lemma2.1},  we obtain $\rho(G_{14})>\rho(G_{13})$.
  Let $G_{15}$ be the hypergraph obtained from $G_{14}$ by applying the edge-releasing operation on $e_1$ at $v_1$ in such a way that $e_1$ becomes a pendent edge.
  By Lemma \ref{lemma2.3}, we have $\rho(G_{15})>\rho(G_{14})$.
  Obviously, $G_{15}$ does not have multiple edges and $G_{15}\in\bar{{\mathcal{U}}}_{1}(n,k,3)$.
  By Corollary \ref{corollary3.1}, we get $\rho(A_{n,k})\geq\rho(G_{15})$ with the equality iff $G_{15}\cong A_{n,k}$.
   By the above proofs, in  Case (i), we obtain $\rho(A_{n,k})>\rho(G)$ for $G\in\bar{\mathcal{U}}_{2,1}(n,k,3)$.  Therefore, we get Lemma \ref{lemma33.4}.

 \noindent\textbf{Case (ii).} $x_{v_{1}}\geq x_{w_{1}}$.

  Let $G_{16}$ be the hypergraph obtained from $G_{13}$ by moving $\dot{g}$ from $w_1$ to $v_1$ in such a way that $G_{16}$ has the perfect matching $M(G)$.
  Obviously, $G_{16}\in \bar{\mathcal{U}}_{2,1}(n,k,3)$.
  By Lemma \ref{lemma2.1}, we obtain $\rho(G_{16})>\rho(G_{13})$.
  If $G_{16}\cong B_{n,k}$, then by Lemma \ref{lemma33.3} and  the above proofs, we get $\rho(A_{n,k})>\rho(B_{n,k})> \rho(G)$ for $G\in\bar{\mathcal{U}}_{2,1}(n,k,3)$. Namely, Lemma \ref{lemma33.4} holds.
  Otherwise, we assume $G_{16}\ncong B_{n,k}$.
  By repeatedly using the same procedure as those in Cases (i) and (ii),
  we finally obtain $\max\{\rho(A_{n,k}),\rho(B_{n,k})\}>\rho(G_{16})$.
  Therefore,  it follows from  $\rho(A_{n,k})>\rho(B_{n,k})$ (by Lemma \ref{lemma33.3}) that $\rho(A_{n,k})> \rho(G_{16})$.
   Thus, in Case (ii), we get $\rho(A_{n,k})>\rho(G)$ for $G\in\bar{\mathcal{U}}_{2,1}(n,k,3)$.

   By combining the above proofs, we have  Lemma \ref{lemma33.4}. ~~$\Box$

%
%

 \begin{lemma}\label{lemma33.5}
  We have $\rho(D_{n,k})>\rho(A_{n,k})$ for  $n\geq 9k(k-1)$ and $k\geq 3$,
 where  $A_{n,k}$ and $D_{n,k}$ are shown in
  Fig.\ \ref{fig-dmr}(b) and  Fig.\ \ref{fig-1}(b), respectively.
\end{lemma}

  \noindent\textbf{Proof}:  We construct a weighted incidence matrix $B_{D_{n,k}}$ for $D_{n,k}$ as follows.
 Let $B_{D_{n,k}}(v,e)=0$ if $v\not\in e$;
   $B_{D_{n,k}}(v,e)=1$ if $v\in e$ and  $v$ is a core vertex;
  $B_{D_{n,k}}(v,e)=\alpha$ if $v\in e$, $e$ is a pendent edge and $d_{D_{n,k}}(v)=2$;
 $B_{D_{n,k}}(v,e)=1-\alpha$ if $v\in e$, $v\neq v_2,v_3$, $e$ is not a pendent edge, and $d_{D_{n,k}}(v)=2$;
 $B_{D_{n,k}}(v_1,g_i)=\dfrac{\alpha}{(1-\alpha)^{k-1}}$ for $i=1,\cdots ,m-2$ and $m\geq 9$;
 $B_{D_{n,k}}(v_1,{e}_{1})=B_{D_{n,k}}(v_1,{e}_{3})=x_{2}$;
 $B_{D_{n,k}}(v_2,{e}_{1})=B_{D_{n,k}}(v_3,{e}_{3})=y_{2}$; and
$B_{D_{n,k}}(v_2,{e}_{2})=B_{D_{n,k}}(v_3,{e}_{2})=c_{1}$,
%
%
 where $y_{2}$, $c_{1}$, $x_{2}$, and $\alpha$ satisfy (\ref{a})--(\ref{c}) as follows:
 \begin{numcases}{}
 y_{2}+c_{1}=1,& \label{a}\\
 2x_{2}+\alpha+\dfrac{(m-2)\alpha}{(1-\alpha)^{k-1}}=1,&\label{cwh}\\
 x_{2}y_{2}(1-\alpha)^{k-2}=\alpha,&\label{b}\\
 c_{1}^{2}=\alpha.&\label{c}
  \end{numcases}
 We can verify that  $\sum_{e:e\in E_{D_{n,k}}(v)}B_{D_{n,k}}(v,e)=1$ for any $v\in V(D_{n,k})$,
  $\prod_{v:v\in e}B_{D_{n,k}}(v,e)=\alpha$ for any $e\in E(D_{n,k})$, and $B_{D_{n,k}}$ is consistent.
   Thus, $D_{n,k}$ is consistently $\alpha$-normal.
    By Lemma \ref{lemma4}(i), we have $\rho({D_{n,k}})=\alpha^{-\frac{1}{k}}$.

  We construct a weighted incidence matrix $B_{A_{n,k}}$ for $A_{n,k}$ as follows.
 Let $B_{A_{n,k}}(v,e)=0$ if $v\not\in e$;
   $B_{A_{n,k}}(v,e)=1$ if $v\in e$ and $v$ is a core vertex;
  $B_{A_{n,k}}(v,e)=\alpha$ if $v\in e$, $e$ is a pendent edge and $d_{A_{n,k}}(v)=2$;
 $B_{A_{n,k}}(v,e)=1-\alpha$ if $v\in e$, $e$ is not a pendent edge and $d_{A_{n,k}}(v)=2$;
 $B_{A_{n,k}}(v_{1},g_i)=\dfrac{\alpha}{(1-\alpha)^{k-1}}$ for $i=1,\cdots ,m-3$ and $m\geq 9$;
 $B_{A_{n,k}}(v_{1},e_{1})=B_{A_{n,k}}(v_{1},e_{3})=x_{3}$;
 $B_{A_{n,k}}(v_{2},e_{1})=B_{A_{n,k}}(v_{3},e_{3})=y_{3}$; and
$B_{A_{n,k}}(v_{2},e_{2})=B_{A_{n,k}}(v_{3},e_{2})=c_{2}$,
%
  where $y_{3}$, $x_{3}$, $c_{2}$, and $\alpha$ satisfy (\ref{d})--(\ref{f}) as follows:
  \begin{numcases}{}
 y_{3}+c_{2}+\alpha=1,& \label{d}\\
 x_{3}y_{3}(1-\alpha)^{k-2}=\alpha,&\label{e}\\
 c_{2}^{2}(1-\alpha)^{k-2}=\alpha.&\label{f}
 \end{numcases}
 We can verify that  $\sum_{e:e\in E_{A_{n,k}}(v)}B_{A_{n,k}}(v,e)=1$ for any $v\in V(A_{n,k})\setminus\{v_1\}$ and
  $\prod_{v:v\in e}B_{A_{n,k}}(v,e)=\alpha$ for any $e\in E(A_{n,k})$.
 Next, we will prove $\sum_{e:e\in E_{A_{n,k}}(v_1)}B_{A_{n,k}}(v_1,e)< 1$.

 Since $\alpha>0$ and $1-\alpha>0$, we have $0<\alpha<1$.
  It follows from (\ref{c}) that $c_{1}=\sqrt{\alpha}$.
  Combining (\ref{a}) and $0<\alpha<1$,  we obtain $y_{2}=1-\sqrt{\alpha}<1-\alpha$.
   From   $y_{2}<1-\alpha$, $y_{2}>0$ and (\ref{b}),  we have $x_{2}=\dfrac{\alpha}{y_{2}(1-\alpha)^{k-2}}>\dfrac{\alpha}{(1-\alpha)^{k-1}}$.
   Substituting $x_{2}>\dfrac{\alpha}{(1-\alpha)^{k-1}}$  into (\ref{cwh}), we get $m\alpha<(1-\alpha)^{k}$.
    It follows from $m\alpha<(1-\alpha)^{k}$ and (\ref{f}) that  $c_{2}<\frac{1-\alpha}{\sqrt{m}}$.
 Thus, for $m\geq9$, we have
 \begin{align}
  &x_{3}-x_{2}=\dfrac{\alpha}{\left(1-\alpha \right)^{k-2}}\left(\frac{1}{y_{3}}-\frac{1}{y_{2}}\right)\label{g}\\
  &=\dfrac{\alpha}{\left(1-\alpha\right)^{k-1}}\left(\frac{1-\alpha}{y_{3}}-\frac{1-\alpha}{y_{2}}\right)\nonumber\\
  &=\dfrac{\alpha}{\left(1-\alpha \right)^{k-1}}\left(\frac{1-\alpha}{1-c_{2}-\alpha}-\frac{1-\alpha}{1-c_{1}}\right)\label{h}\\
  &<\dfrac{\alpha}{\left(1-\alpha \right)^{k-1}}\left(\frac{\sqrt{m}}{\sqrt{m}-1}-\left(1+\sqrt{\alpha}\right)\right)\label{i}\\
  &<\dfrac{\alpha}{\left(1-\alpha \right)^{k-1}}\left(\frac{\sqrt{m}}{\sqrt{m}-1}-1\right)\nonumber\\
  &=\dfrac{\alpha}{\left(1-\alpha \right)^{k-1}}\times\frac{1}{\sqrt{m}-1}\nonumber\\
  &\leq\dfrac{\alpha}{2\left(1-\alpha \right)^{k-1}}\label{j},
  \end{align}
 where (\ref{g}) follows from (\ref{b}) and  (\ref{e}),
  (\ref{h}) is  deduced  from  (\ref{a}) and  (\ref{d}), and
  (\ref{i}) is obtained from  $c_{2}<\frac{1-\alpha}{\sqrt{m}}$ and $c_{1}=\sqrt{\alpha}$.
  Therefore, it follows from (\ref{g})--(\ref{j}) that  $x_{3}<x_{2}+\dfrac{\alpha}{2(1-\alpha)^{k-1}}$ for $m\geq 9$.
  Thus,  we obtain

 \begin{align}
  &\sum_{e:e\in E_{A_{n,k}}(v_1)}B_{A_{n,k}}(v_1,e)\nonumber \\
  &=2x_{3}+\alpha+\dfrac{(m-3)\alpha}{(1-\alpha)^{k-1}}\nonumber\\
  &<2x_{2}+\dfrac{\alpha}{(1-\alpha)^{k-1}}+\alpha+\dfrac{(m-3)\alpha}{(1-\alpha)^{k-1}}\nonumber\\
  &=2x_{2}+\alpha+\dfrac{(m-2)\alpha}{(1-\alpha)^{k-1}}\nonumber\\
  &=1,\label{iwh}
     \end{align}
     where (\ref{iwh}) follows from (\ref{cwh}).
 Thus, for $m\geq9$, $A_{n,k}$ is strictly $\alpha$-subnormal.
 Therefore, by Lemma \ref{lemma4}(ii), we have $\rho({A_{n,k}})<\alpha^{-\frac{1}{k}}$.

 In conclusion, it follows from Lemma \ref{lemma4} that $\rho(D_{n,k})>\rho(A_{n,k})$, where $n\geq 9k(k-1)$ and $k\geq 3$.
  ~~$\Box$

  It should be noted that  the methods proposed in Lemma \ref{lemma33.5} can not be used to compare the relationship between  $\rho(D_{n,k})$ and $\rho(A_{n,k})$ when  $n=mk(k-1)$ with $3\leq m\leq8$ since (\ref{j}) holds for  $m\geq 9$.

   By the methods similar to those for  Lemma \ref{lemma33.1},  we have Lemma  \ref{lemma33.2} as follows.

  \begin{lemma}\label{lemma33.2}
  Let $G\in\bar{\mathcal{U}}_{2,2}(n,k,3)$, where $n\geq 2k(k-1)$ and $k\geq 3$.
   We have $\rho(D_{n,k})\geq\rho(G)$  with the equality iff $G\cong D_{n,k}$.
 \end{lemma}


 \begin{corollary}\label{corollary3.2}
 Let $G \in \mathcal{U}_{2}(n,k,l)$, where $n=mk(k-1)$ and $m,k, l\geq 3$.
(i). If $3\leq m\leq8$,   we have $\max\{\rho(A_{n,k}),\rho(D_{n,k})\}\geq\rho(G)$.
    (ii). If $m\geq 9$, we have $\rho(D_{n,k})\geq\rho(G)$ with the equality iff $G\cong D_{n,k}$.
 \end{corollary}

  \noindent\textbf{Proof}: Let $n=mk(k-1)$ and $m,k, l\geq 3$.
   By  Corollary \ref{corollary3.1}, we have $\rho(A_{n,k})\geq\rho(G)$ with the equality iff $G\cong A_{n,k}$, where  $G\in \mathcal{U}_{1}(n,k,l)$. By Lemmas \ref{lemma33.4} and  \ref{lemma33.2}, we have $\max\{\rho(A_{n,k}),\rho(D_{n,k})\}\geq\rho(G)$ for $G\in\bar{\mathcal{U}}_{2}(n,k,3)$  since
  \textcolor{black}{$\bar{\mathcal{U}}_{2,1}(n,k,3)\cup \bar{\mathcal{U}}_{2,2}(n,k,3)=\bar{\mathcal{U}}_{2}(n,k,3)$.}
   Furthermore,
  \textcolor{black}{by Lemmas \ref{lemma3.1}, \ref{lemma3.3} and Corollary \ref{corollary3.1},}
  we obtain $\max\{\rho(A_{n,k}),\rho(D_{n,k})\}\geq\rho(G)$ for $G\in\mathcal{U}_{2}(n,k,l)$. Thus, when $3\leq m\leq8$, we get Corollary \ref{corollary3.2} (i).
   For $m\geq 9$, by Lemma \ref{lemma33.5}, we have  $\rho(D_{n,k})>\rho(A_{n,k})$ for  $n\geq 9k(k-1)$. Therefore, we get Corollary \ref{corollary3.2} (ii).
       ~~$\Box$

  By
  \textcolor{black}{Corollaries \ref{corollary3.1} and \ref{corollary3.2},}
  we obtain Theorem \ref{theorem3.6} as follows.

  \begin{theorem}\label{theorem3.6}
  Let $G\in \mathcal{U}(n,k)$, where  $n=mk(k-1)$, $m\geq2$ and $k\geq 3$.
  (i). If $m=2$, we have $G\cong B_{n,k}\cong D_{n,k}$ and $\rho(G)=\rho(B_{n,k})=\rho(D_{n,k})$.
   (ii). If $3\leq m\leq8$,    $\max\{\rho(A_{n,k}),\rho(D_{n,k})\}\geq\rho(G)$.
    (iii). If $m\geq 9$, $\rho(D_{n,k})\geq\rho(G)$ with the equality iff $G\cong D_{n,k}$.
 \end{theorem}

 \section{The hypergraph with the maximal spectral radius among $\Gamma(n,k)$}

 In  Section 4, we will get the hypergraph with the maximal spectral radius among $ \Gamma(n,k)$, where $n=mk(k-1)$, $m\geq 1$ and $k\geq 3$.
 Some necessary definitions are given as follows.

 Let $G \in \Gamma(n,k)$. The unique cycle in $G$ has  two edges which share two common vertices.
 We denote the cycle in $G$ by $C_{2}$ and the two edges contained in $C_{2}$ by $\widetilde{e}_{1}$ and $\widetilde{e}_{2}$,
 where $\widetilde{e}_{1}=\{u_{1},u_{1,1},...,u_{1,k-2},u_{2}\}$ and $\widetilde{e}_{2}=\{u_{1},u_{2,1},...,u_{2,k-2},u_{2}\}$.
 $C_{2}$ is shown in Fig.\ \ref{fig-a}(a).
  According to the fact whether $C_2$ of $G $ has  one perfect matching edge or not, we classify  $\Gamma(n,k)$ into two types:
 (1)  $C_2$ has one perfect matching edge, and (2) $\widetilde{e}_{1}$ and  $\widetilde{e}_{2}$ in $C_2$  are not perfect matching edges.
 The hypergraphs  in   $\Gamma(n,k)$ can be divided into two subsets according to   Types (1) and (2).
    We denote
 $\Gamma(n,k)=\Gamma_{1}(n,k)\cup \Gamma_{2}(n,k)$,
 where  all the hypergraphs in  $\Gamma_{1}(n,k)$  and $\Gamma_{2}(n,k)$ have  Types (1) and (2), respectively.
 Obviously,  for all the hypergraphs in $\Gamma_{1}(n,k)$  and $\Gamma_{2}(n,k)$, we have $m\geq 1$ and  $m\geq 2$, respectively.




 Let $\bar{\Gamma}_{1}(n,k)$ be a subset of $\Gamma_{1}(n,k)$ in which each hypergraph satisfies  three conditions:
 (i) each vertex in $\widetilde{e}_{1}\setminus\{u_{1},u_{2}\}$ of $C_2$ is a core vertex;
(ii) each vertex in $\widetilde{e}_{2}\setminus\{u_{1},u_{2}\}$ of  $C_2$ must be attached by a pendent edge; and
 (iii) at most one of the vertices  (denoted by $v$) of $\widetilde{e}_{2}$ of  $C_2$  may or may not be attached
  by a hypertree which has at least $k\geq 3$ edges.
 We further classify  $\bar{\Gamma}_{1}(n,k)$  into two subsets which are denoted by $\bar{\Gamma}_{1,1}(n,k)$ and $\bar{\Gamma}_{1,2}(n,k)$, where the hypergraphs in $\bar{\Gamma}_{1,1}(n,k)$ satisfy  that $v=u_{1}$ or $v=u_{2}$ and
 the hypergraphs in $\bar{\Gamma}_{1,2}(n,k)$ satisfy that  $v$ is one of the vertices in
  $\{u_{2,1},\cdots,u_{2,k-2}\}$.

  Let $\bar{\Gamma}_{2}(n,k)$ be a subset of
  \textcolor{black}{$\Gamma_{2}(n,k)$}
  in which each hypergraph satisfies  two conditons:
   (i) each vertex in $C_2$ must be attached by a pendent edge, and
  (ii) at most one of the vertices in $\widetilde{e}_{1}\cap\widetilde{e}_{2}=\{u_1, u_2\}$ of  $C_2$  is attached by a hypertree which has at least $k$ edges, where $k\geq 3$.

\begin{figure}
\thicklines
\begin{tikzpicture}
\draw[thick,black,rotate=45](1.85,-0.7) ellipse(1.3 and 0.2); 
\draw[thick,black,rotate=135](0.4,-0.75) ellipse(1.3 and 0.2); 
\draw[thick,blue](-1.3,1.5) ellipse(1.0 and 0.15); 
\draw[black](-0.5,1.5) circle(0.02);  
\draw[thick,blue](3.35,1.5)  ellipse(1.0 and 0.15); 
\draw[black](2.55,1.5) circle(0.02);  
\draw[thick,blue](-0.5,0.55) ellipse(1.0 and 0.15); 
\draw[black](0.35,0.55) circle(0.02);  
\draw[thick,blue](2.5,0.55)  ellipse(1.0 and 0.15); 
\draw[black](1.7,0.55) circle(0.02);  
\draw[thick,black](0.8,0) ellipse(1.5 and 0.2); 
\draw[thick,red](3,0)  ellipse(1.5 and 0.2); 
\fill(0.3,1.8) circle(0.5pt);
\fill(0.9,1.8) circle(0.5pt);
\fill(1.5, 1.8) circle(0.5pt); 
\fill(-0.25, 1.3) circle(0.5pt);
\fill(0.15, 0.9) circle(0.5pt);
\fill(-0.05, 1.1) circle(0.5pt);
\fill(2.05, 1) circle(0.5pt);
\fill(1.85, 0.8) circle(0.5pt);
\fill(2.25, 1.2) circle(0.5pt);
\draw[thick,blue,rotate=270](0.9,-0.35) ellipse(1 and 0.15); 
\draw[thick,blue,rotate=270](0.9,1) ellipse(1 and 0.15); 
\draw[black](-0.35,0) circle(0.02);  
\fill(0.40,0) circle(0.02);  
\fill(0.01,0) circle(0.02);  
\fill(0.75,0) circle(0.02); 
\fill(0.40,-1) circle(0.02);  
\fill(0.01,-1) circle(0.02);  
\fill(0.75,-1) circle(0.02); 
\draw[black](1.0,0)  circle(0.02);  
\draw[black](1.7,0)  circle(0.02);  
\draw[black](2.05,0)  circle(0.02);  
\fill(2.7,0) circle(0.02);   
\fill(3.3,0) circle(0.02);   
\fill(4,0) circle(0.02);   
\draw(1.05,2.5) node{{\scriptsize $m-1$}};
 \draw(1.,2) node{$\overbrace{~~~~~~~~~~~~~~~~~~~}$};
\draw(-1,2) node{{\scriptsize $g_1$}};
\draw(3,2) node{{\scriptsize $g_{m-1}$}};
\draw(-1,0) node{{\scriptsize $\widetilde{e}_{2}$}};
\draw(4.75,0) node{{\scriptsize $\widetilde{e}_{1}$}};
\draw(2.2,-0.4) node{{\scriptsize $u_1$}};
\draw(1.57,-0.4) node{{\scriptsize $u_2$}};
\draw(0.6,-0.4) node{{\scriptsize $u_{2,1}$}};
\draw(2,-3) node{{\scriptsize $(a)$}};
\end{tikzpicture}
\begin{tikzpicture}
\draw[thick,black,rotate=45](2.4,-1.2) ellipse(1.3 and 0.2); 
\draw[thick,black,rotate=135](-0.05,-1.2) ellipse(1.3 and 0.2); 
\draw[thick,blue](-0.42,1.2) ellipse(1.0 and 0.15); 
\draw[black](0.42,1.2) circle(0.02);  
\draw[thick,blue](3.7,1.2)  ellipse(1.0 and 0.15); 
\draw[black](2.9,1.2) circle(0.02);  
\draw[thick,blue](0.19,0.55) ellipse(1.0 and 0.15); 
\draw[black](1.1,0.55) circle(0.02);  
\draw[thick,blue](3.1,0.55)  ellipse(1.0 and 0.15); 
\draw[black](2.22,0.55) circle(0.02);  
\draw[thick,black](0.8,0) ellipse(1.5 and 0.2); 
\draw[thick,black](3,0)  ellipse(1.5 and 0.2); 
\fill(1.0,1.8) circle(0.5pt);
\fill(1.5,1.8) circle(0.5pt);
\fill(2.0, 1.8) circle(0.5pt); 
\fill(0.8, 1) circle(0.5pt);
\fill(0.88, 0.9) circle(0.5pt);
\fill(0.98, 0.8) circle(0.5pt);
\fill(2.4, 0.8) circle(0.5pt);
\fill(2.5, 0.9) circle(0.5pt);
\fill(2.6, 1.0) circle(0.5pt);
\draw[thick,blue,rotate=270](0.9,-0.35) ellipse(1 and 0.15); 
\draw[thick,blue,rotate=270](0.9,1) ellipse(1 and 0.15); 
\draw[thick,blue,rotate=270](0.9,1.7) ellipse(1 and 0.15); 
\draw[thick,blue,rotate=270](0.9,2.05) ellipse(1 and 0.15); 
\draw[thick,blue,rotate=270](0.9,2.8) ellipse(1 and 0.15); 
\draw[thick,blue,rotate=270](0.9,4.15) ellipse(1 and 0.15); 
\draw[black](-0.35,0) circle(0.02);  
\fill(0.40,0) circle(0.02);  
\fill(0.01,0) circle(0.02);  
\fill(0.75,0) circle(0.02); 
\fill(0.40,-0.7) circle(0.02);  
\fill(0.01,-0.7) circle(0.02);  
\fill(0.75,-0.7) circle(0.02); 
\draw[black](1.0,0)  circle(0.02);  
\draw[black](1.7,0)  circle(0.02);  
\draw[black](2.05,0)  circle(0.02);  
\fill(3.2,0) circle(0.02);   
\fill(3.6,0) circle(0.02);   
\fill(4,0) circle(0.02);   
\fill(3.1,-0.7) circle(0.02);   
\fill(3.5,-0.7) circle(0.02);   
\fill(3.85,-0.7) circle(0.02);   
\draw[black](4.15,-0.04) circle(0.02);   
\draw[black](2.8,0) circle(0.02);  
\draw(-0.2,2) node{{\scriptsize $g_1$}};
\draw(3.7,2) node{{\scriptsize $g_{m-2}$}};
\draw(-1,0) node{{\scriptsize $\widetilde{e}_{1}$}};
\draw(4.75,0) node{{\scriptsize $\widetilde{e}_{2}$}};
\draw(2.37,-0.4) node{{\scriptsize $u_2$}};
\draw(1.4,-0.4) node{{\scriptsize $u_1$}};
\draw(1.5,2.5) node{{\scriptsize $m-2$}};
 \draw(1.63,2) node{$\overbrace{~~~~~~~~~~~~~~~~~~~}$};
 \draw(1.9,-3) node{{\scriptsize $(b)$}};
\end{tikzpicture}
\caption{\label{fig-5} (a) $J_{n,k}$ and (b) $L_{n,k}$}
\end{figure}
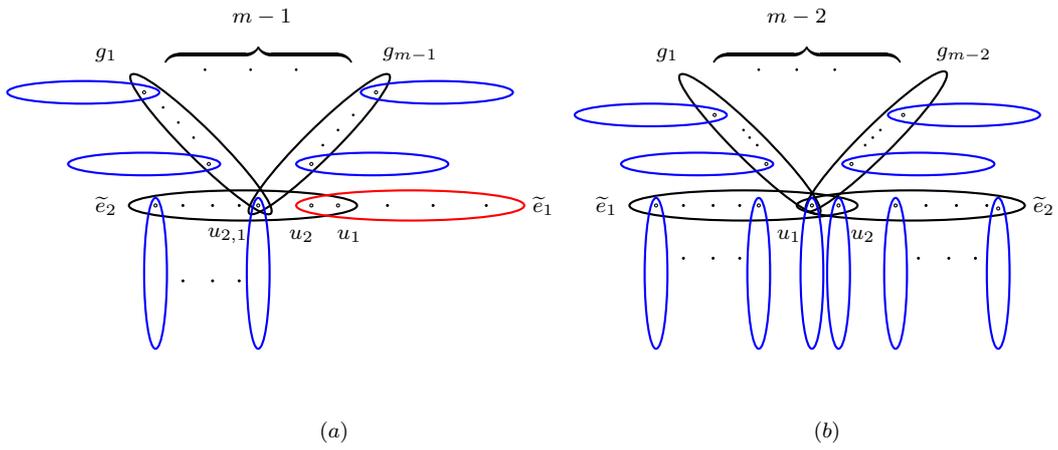

 Let $I_{n,k}$ be  the hypergraph obtained from $C_{2}(u_{1},u_{0})S_{m-1,k}$ by attaching one pendent edge  at each vertex of $C_{2}(u_{1},u_{0})S_{m-1,k}$ (except for all the vertices  of $\widetilde{e}_{1}$), where  $m\geq 1$ and $k\geq3$.
    Let  $J_{n,k}$ be the hypergraph obtained from $C_{2}(u_{2,1},u_{0})S_{m-1,k}$ by attaching one pendent edge  at each vertex of $C_{2}(u_{2,1},u_{0})S_{m-1,k}$ (except for all the vertices  of $\widetilde{e}_{1}$), where  $m\geq 1$ and $k\geq3$.
   Obviously, when $n=k(k-1)$, $I_{n,k}\cong J_{n,k}$.
    $I_{n,k}$ and $J_{n,k}$ are  shown in Fig.\ \ref{fig-a}(b) and Fig.\ \ref{fig-5}(a), respectively.
   Let $L_{n,k}$ be the hypergraph obtained from $C_{2}(u_{1},u_{0})S_{m-2,k}$ by attaching one pendent edge at each vertex of $C_{2}(u_{1},u_{0})S_{m-2,k}$, where  $m\geq2$ and $k\geq3$. $L_{n,k}$ is shown in Fig.\ \ref{fig-5}(b).
  Obviously,  $I_{n,k}\in\bar{\Gamma}_{1,1}(n,k)$, $J_{n,k}\in\bar{\Gamma}_{1,2}(n,k)$, and $L_{n,k}\in\bar{\Gamma}_{2}(n,k)$.

 To obtain the hypergraph with the maximal spectral radius in $\Gamma(n,k)=\Gamma_{1}(n,k) \cup \Gamma_{2}(n,k)$ (as shown in  Theorem \ref{theorem4.6}), we introduce several lemmas  first.
  We propose Lemmas \ref{lemmawang1}--\ref{lemma44.3} to get the hypergraph with the maximal spectral radius in $\Gamma_{1}(n,k)$ (as shown in  Corollary \ref{corollary4.2}).
   Lemmas \ref{lemma4.4} and \ref{lemma44.1} are deduced to obtain
    the  hypergraph with the maximal spectral radius in $\Gamma_{2}(n,k)$ (as shown in  Corollary \ref{corollary4.1}).

 By the methods similar to those for Lemma \ref{lemmawang}, we have Lemma \ref{lemmawang1} as follows.

\begin{lemma}\label{lemmawang1}
  Let $G\in \Gamma_i(n,k)$, where $i=1,2$, $n\geq k(k-1)$ and $k\geq3$.
  Let $e$ be a perfect matching edge of $G$, $e$ is not an edge on $C_2$ of $G$ and $e$ is not a pendent edge.
  Let $G'_0$ be the hypergraph obtained from $G$ by applying the edge-releasing operation on $e$ at a vertex of $e$ such that $e$ of $G'_0$ is a pendent edge.
 We have $\rho(G'_0)>\rho(G)$, where $G'_0\in \Gamma_i(n,k)$ and $i=1,2$.
 \end{lemma}

 \begin{lemma}\label{lemma4.3}
  Let $G\in\Gamma_{1}(n,k)$, where  $n\geq k(k-1)$ and $k\geq 3$.
  There exists one hypergraph $\ddot{G}\in\bar{\Gamma}_{1}(n,k)$ such that $\rho(\ddot{G})\geq\rho(G)$ with the equality iff $G\cong \ddot{G}$.
 \end{lemma}

 \noindent\textbf{Proof}:  Let $n\geq k(k-1)$ and $k \geq 3$.
  Let $G\in\Gamma_{1}(n,k)$.
  According to the definition  of $\Gamma_{1}(n,k)$, we suppose that $\widetilde{e}_{1}$ of $G$ is a perfect matching edge.
  In $\widetilde{e}_{1}\setminus \{u_1,u_2\}$, if there exists one vertex having degree not less than 2,  without loss of generality, we suppose that this vertex is $u_{1,1}$.
  By the methods similar to those for the first paragraph in Lemma \ref{lemma3.3}, we finally obtain a hypergraph (denoted by $G_{17}$) in $\Gamma_{1}(n,k)$ satisfying $\rho(G_{17})\geq\rho(G)$ with the equality iff $G\cong G_{17}$
  and each vertex in $\widetilde{e}_{1}\backslash\{u_{1},u_{2}\}$ of $G_{17}$ is a core vertex.

  If $G_{17}\in\bar{\Gamma}_{1}(n,k)$, then we get Lemma \ref{lemma4.3}.
  Otherwise, by Lemma \ref{lemmawang1},
       we  obtain a hypergraph $\tilde{G}$
   satisfying that
   \textcolor{black}{$\rho(\tilde{G})\geq\rho(G_{17})$ with the equality iff $G_{17}\cong\tilde{G}$,}
   where $\tilde{G}$ satisfies three conditions:
   (i) each vertex in $\widetilde{e}_{1}\setminus\{u_{1},u_{2}\}$ of $C_2$ is a core vertex;
(ii) each vertex in $\widetilde{e}_{2}\setminus\{u_{1},u_{2}\}$ of  $C_2$ must be attached by a pendent edge; and
(iii) the vertex in $\widetilde{e}_{2}$ of  $C_2$  may or may not be attached by a hypertree which has at least $k\geq 3$ edges.

  If $\tilde{G} \in\bar{\Gamma}_{1}(n,k)$, then we get Lemma \ref{lemma4.3}.
   Otherwise, we assume $\tilde{G} \not\in\bar{\Gamma}_{1}(n,k)$. By the definition  of $\bar{\Gamma}_{1}(n,k)$,
   there exist two vertices (denoted by $\tilde{u}_1$ and $\tilde{u}_2$) in $\widetilde{e}_{2}$ of $\tilde{G}$ which have degrees not less than 3. Namely,
 $\tilde{u}_1$ and $\tilde{u}_2$ are attached by hypertrees which have  at least $k\geq 3$  edges.
 Let $f_{1}^{1},\cdots,f_{d_{G}(\tilde{u}_1)-2}^{1}$ and $f_{1}^{2},\cdots,f_{d_{G}(\tilde{u}_2)-2}^{2}$ be all the edges
 which are incident with
 $\tilde{u}_1$ and $\tilde{u}_2$ respectively.
 It is noted that $f_{1}^{1},\cdots,f_{d_{G}(\tilde{u}_1)-2}^{1}$ and $f_{1}^{2},\cdots,f_{d_{G}(\tilde{u}_2)-2}^{2}$ are not perfect matching  edges and all of them do not contain $\widetilde{e}_{2}$.

 Let $\bm{x}$ be the principal eigenvector of $\tilde{G}$ corresponding to $\rho(\tilde{G})$. Without loss of generality, we suppose $x_{\tilde{u}_1}\geq x_{\tilde{u}_2}$.
 Let $G_{18}$ be the hypergraph obtained from $\tilde{G}$ by removing $f_{1}^{2},\cdots,f_{d_{G}(\tilde{u}_2)-2}^{2}$ from
 $\tilde{u}_2$ to $\tilde{u}_1$. 
  By Lemma \ref{lemma2.1}, we obtain $\rho(G_{18})>\rho(\tilde{G})$.
  By repeatedly using the same procedure as above,
  we can find a hypergraph $\ddot{G}\in\bar{\Gamma}_{1}(n,k)$ such that $\rho(\ddot{G})\geq\rho(G_{18})>\rho(G)$ with the equality iff $\ddot{G}\cong G_{18}$.
 Therefore, we get Lemma \ref{lemma4.3}.~~$\Box$

  To  obtain the hypergraph with the maximal spectral radius in $\bar{\Gamma}_{1,1}(n,k)$ with  $n\geq 2k(k-1)$ and $k\geq3$
  (as shown in Lemma \ref{lemma44.5}), we introduce Lemmas \ref{lemma444} and \ref{lemma44.4} first.

  \begin{lemma}\label{lemma444}
  We have $\rho(L_{n,k})>\rho(A_{n,k})$ and  $\rho(L_{n,k})>\rho(D_{n,k})$, where $n\geq 2k(k-1)$ and $k\geq3$.
 \end{lemma}

  \noindent\textbf{Proof}: Let $n\geq 2k(k-1)$ and $k\geq3$. Let $\bm{x}$ be the principal eigenvector corresponding to $\rho(A_{n,k})$.
  In $A_{n,k}$, if $x_{v_1}\geq x_{v_3}$, then let
  $H_{1}$ be the hypergraph obtained from $A_{n,k}$ by removing the edge $e_{2}$ of $A_{n,k}$ from $v_3$ to $v_1$.
  By Lemma \ref{lemma2.1},  we have $\rho(H_{1})>\rho(A_{n,k})$.
   In $A_{n,k}$, if $x_{v_3}> x_{v_1}$, then let $H_{2}$ be the hypergraph obtained from $A_{n,k}$ by removing all the edges of $e_{1},g_{1},\cdots,g_{m-3} (m\geq 3)$ from $v_1$ to $v_3$.
  By Lemma \ref{lemma2.1}, we obtain $\rho(H_{2})>\rho(A_{n,k})$.
  Obviously, $H_{1}\cong H_{2}\cong L_{n,k}$.
   Therefore, we get $\rho(L_{n,k})>\rho(A_{n,k})$ for $n\geq 3k(k-1)$ and $k\geq3$.

  Let $\bm{y}$ be the principal eigenvector corresponding to $\rho(D_{n,k})$.
   By the symmetry of the vertices of $D_{n,k}$, we have $y_{v_2}=y_{v_3}$.
   Let  $H_{3}$ be the hypergraph obtained from $D_{n,k}$  by removing
   ${e}_{3}$ from $v_3$ to $v_2$. Obviously, $H_{3}\cong L_{n,k}$.
    Therefore, by Lemma \ref{lemma2.1},  we get $\rho(L_{n,k})>\rho(D_{n,k})$ for $n\geq 2k(k-1)$ and $k\geq3$.
  ~~$\Box$

 \begin{lemma}\label{lemma44.4}  
 We have $\rho(L_{n,k})> \rho(I_{n,k})$ for $n\geq 2k(k-1)$ and $k\geq 3$, where $I_{n,k}$ and $L_{n,k}$   are shown in Fig.\ \ref{fig-a}(b) and Fig.\ \ref{fig-5}(b), respectively.
\end{lemma}

 \noindent\textbf{Proof}: Let $n\geq 2k(k-1)$ and $k\geq 3$. We construct a weighted incidence matrix $B_{L_{n,k}}$ for $L_{n,k}$ as follows.
 Let $B_{L_{n,k}}(v,e)=0$ if $v\not\in e$;
   $B_{L_{n,k}}(v,e)=1$ if $v\in e$ and $v$ is a core vertex;
  $B_{L_{n,k}}(v,e)=\alpha$ if $v\in e$, $e$ is a pendent edge and $d_{L_{n,k}}(v)=2$;
 $B_{L_{n,k}}(v,e)=1-\alpha$ if $v\in e$, $e$ is not a pendent edge and $d_{L_{n,k}}(v)=2$;
 $B_{L_{n,k}}(u_1,g_i)=\dfrac{\alpha}{(1-\alpha)^{k-1}}$ for $i=1,\cdots ,m-2$ and $m\geq 2$;
 $B_{L_{n,k}}(u_1,\widetilde{e}_{1})=B_{L_{n,k}}(u_1,\widetilde{e}_{2})=x_{4}$; and
 $B_{L_{n,k}}(u_2,\widetilde{e}_{1})=B_{L_{n,k}}(u_2,\widetilde{e}_{2})=y_{4}$,
 where $x_{4}$, $y_{4}$ and $\alpha$ satisfy (\ref{2.1})--(\ref{2.3}) as follows:
 \begin{numcases}{}
 2x_{4}+\alpha   +\dfrac{(m-2)\alpha   }{(1-\alpha   )^{k-1}}=1,&\label{2.1} \\
 2y_{4}+\alpha   =1,&\label{2.2}\\
 x_{4}y_{4}(1-\alpha )^{k-2}=\alpha   .&\label{2.3}
 \end{numcases}
 We can check that  $\sum_{e:e\in E_{L_{n,k}}(v)}B_{L_{n,k}}(v,e)=1$ for any $v\in V(L_{n,k})$,
  $\prod_{v:v\in e}B_{L_{n,k}}(v,e)=\alpha$ for any $e\in E(L_{n,k})$, and $B_{L_{n,k}}$ is consistent.
   Thus, $L_{n,k}$ is consistently $\alpha$-normal.
    By Lemma \ref{lemma4}(i), we have $\rho({L_{n,k}})=\alpha^{-\frac{1}{k}}$.


  We construct a weighted incidence matrix $B_{I_{n,k}}$ for $I_{n,k}$ as follows.
 Let $B_{I_{n,k}}(v,e)=0$ if $v\not\in e$;
   $B_{I_{n,k}}(v,e)=1$ if $v\in e$  and $v$ is a core vertex;
  $B_{I_{n,k}}(v,e)=\alpha$ if $v\in e$, $e$ is a pendent edge and $d_{I_{n,k}}(v)=2$;
 $B_{I_{n,k}}(v,e)=1-\alpha$ if $v\in e$ ($v\neq u_2$), $e$ is not a pendent edge and $d_{I_{n,k}}(v)=2$;
 $B_{I_{n,k}}(u_1,g_i)=\dfrac{\alpha}{(1-\alpha)^{k-1}}$ for $i=1,\cdots ,m-1$ and $m\geq 2$;
 $B_{I_{n,k}}(u_1,\widetilde{e}_{1})=x_{5}$;
 $B_{I_{n,k}}(u_2,\widetilde{e}_{1})=y_{5}$;
 $B_{I_{n,k}}(u_1,\widetilde{e}_{2})=x_{4}$; and
 $B_{I_{n,k}}(u_2,\widetilde{e}_{2})=y_{4}$,
  where $x_{4}$, $y_{4}$, $x_{5}$, $y_{5}$, and $\alpha $ satisfy (\ref{2.4}) and (\ref{2.5}) as follows:
  \begin{numcases}{}
 x_{4}+x_{5}+\dfrac{(m-1)\alpha  }{(1-\alpha  )^{k-1}}=1,&\label{2.4} \\
 y_{4}+y_{5}=1.&\label{2.5}
 \end{numcases}
 We can check that  $\sum_{e:e\in E_{I_{n,k}}(v)}B_{I_{n,k}}(v,e)=1$ for any $v\in V(I_{n,k})$ and
  $\prod_{v:v\in e}B_{I_{n,k}}(v,e)=\alpha$ for any $e\in E(I_{n,k})$ and $e\neq \widetilde{e}_{1}$.
  Next, we prove $\prod_{v: v\in \widetilde{e}_{1}}B_{I_{n,k}}(v,\widetilde{e}_{1})>\alpha$.
 We have
\begin{align}
  &\prod_{v:v\in \widetilde{e}_{1}}B_{I_{n,k}}(v,\widetilde{e}_{1})-\alpha   \nonumber\\
  &=x_{5}y_{5}-\alpha   \nonumber\\
   &=[x_{4}+\alpha   -\dfrac{\alpha   }{(1-\alpha   )^{k-1}}](y_{4}+\alpha   )-\alpha   \label{111}\nonumber\\
   \end{align}
 \begin{align}
  &=x_{4}y_{4}+x_{4}\alpha   +y_{4}\alpha   +\alpha   ^2-\dfrac{\alpha   }{(1-\alpha   )^{k-1}}y_{4}-\dfrac{\alpha   ^2}{(1-\alpha   )^{k-1}}-\alpha   \nonumber\\
  &=\dfrac{\alpha   }{(1-\alpha   )^{k-2}}+\dfrac{2\alpha   ^2}{(1-\alpha   )^{k-1}}+\dfrac{\alpha   -\alpha   ^2}{2}+\alpha   ^2-\dfrac{\alpha   }{2(1-\alpha   )^{k-2}}-\dfrac{\alpha   ^2}{(1-\alpha   )^{k-1}}-\alpha   \label{222}\\
  &=\dfrac{\alpha   }{2(1-\alpha   )^{k-2}}+\dfrac{\alpha   ^2}{(1-\alpha   )^{k-1}}+\dfrac{\alpha   ^2}{2}-\dfrac{\alpha   }{2}\nonumber\\
  &>\dfrac{\alpha   }{2}+\dfrac{\alpha   ^2}{(1-\alpha   )^{k-1}}+\dfrac{\alpha   ^2}{2}-\dfrac{\alpha   }{2}\label{333}\\
  &>\frac{3\alpha   ^2}{2}\label{666}\\
  &>0.\label{676}
     \end{align}
 It is noted that (\ref{111}) follows from $x_{5}=x_{4}+\alpha   -\dfrac{\alpha   }{(1-\alpha   )^{k-1}}$
 (by (\ref{2.1}) and (\ref{2.4})) and $y_{5}=y_{4}+\alpha $
 (by (\ref{2.2}) and (\ref{2.5})).
    Since $y_{4}=\dfrac{1-\alpha   }{2}$ (by (\ref{2.2})) and $x_{4}=\dfrac{2\alpha   }{(1-\alpha   )^{k-1}}$ (by (\ref{2.3}) and $y_{4}=\dfrac{1-\alpha   }{2}$),
    we get (\ref{222}). Since  $0<1-\alpha   <1$, we obtain (\ref{333})--(\ref{676}). Thus, we get $\prod_{v: v\in \widetilde{e}_{1}}B_{I_{n,k}}(v,\widetilde{e}_{1})>\alpha$.

 By the above proofs, we get that   $I_{n,k}$ is strictly $\alpha$-subnormal.
 Therefore, by Lemma \ref{lemma4}(ii), we have $\rho({I_{n,k}})<\alpha^{-\frac{1}{k}}$.
 Thus, by Lemma \ref{lemma4},  we obtain $\rho(L_{n,k})> \rho(I_{n,k})$ for $n\geq 2k(k-1)$ and $k\geq 3$.
 ~~$\Box$

   \begin{lemma}\label{lemma44.5} 
 Let $G\in\bar{\Gamma}_{1,1}(n,k)$, where  $n\geq 2k(k-1)$ with $k\geq3$.
  We have  $\rho(L_{n,k})>\rho(G)$.
 \end{lemma}

 \noindent\textbf{Proof}:
Let $n =mk(k-1)$ with $m\geq 2$ and $k\geq 3$.
 Let $G\in\bar{\Gamma}_{1,1}(n,k)$.
 When $m=2$,  we get Lemma \ref{lemma44.5} since  $G\cong I_{n,k}$ and $\rho(L_{n,k})>\rho(I_{n,k})$ (by Lemma \ref{lemma44.4}).
  Next,  let  $m \geq 3$.

    By the definition of $\bar{\Gamma}_{1,1}(n,k)$, each vertex in $\widetilde{e}_{2}\setminus\{u_{1},u_{2}\}$ of $G$ is attached by a pendent edge and is not attached by  a hypertree  which has at least $k$ ($k\geq 3$) edges, and we assume $u_{1}$ of $C_2$ of $G$ is attached by a hypertree (denoted by $\overline{T}$) which has at least two edges belonging to  $Q(G)$ and $\overline{T}$ maybe has perfect matching edges which are not pendent edges.
  Let  $M(G)$ be the perfect matching of $G$.
    By Lemma \ref{lemmawang1},
   we  get a hypergraph $G_{19}$ such that $\rho(G_{19})\geq\rho(G)$ with the equality iff
 $G\cong G_{19}$, where $G_{19}\in\bar{\Gamma}_{1,1}(n,k)$
 and all the  perfect matching edges of $G_{19}$
 are pendent edges except for $\widetilde{e}_1$ contained in $C_2$ of $G_{19}$.

   If $G_{19}\cong I_{n,k}$,
   then by  Lemma \ref{lemma44.4}, we get  $\rho(L_{n,k})>\rho(I_{n,k})\geq \rho(G)$ and Lemma \ref{lemma44.5} holds.
   Otherwise, we assume $G_{19}\ncong I_{n,k}$.
  It is noted that $u_{1}$ of $C_2$ of $G_{19}$
  is attached by a hypertree (denoted by $\widetilde{T}$) which has at least two edges  belonging to $Q(G_{19})$.
   As $G_{19}\ncong I_{n,k}$, in  $G_{19}$,
   there exists one edge (denoted by $g'=\{z_{1},...,z_{k}\}$) which satisfies the following conditions:
  (i) $g'\in E(\widetilde{T})$;
  (ii) $g'$ is not a pendent edge;
  (iii) $u_1$ and $z_1$ are incident with a common edge (denoted by $g''$);
  and (iv) $g'$ is not adjacent to $u_1$.
    Let $\bm{x}$ be the principal eigenvector of $G_{19}$ corresponding to $\rho(G_{19})$.
  Two cases are considered as follows.
%

\noindent\textbf{ Case (i).} $x_{u_{1}}<x_{z_{1}}$.

  Let $G_{20}$ be the hypergraph obtained from $G_{19}$
  by moving all the edges which are incident with $u_1$
  (except for $\widetilde{e}_1$ and $g''$ of $G_{19}$) from $u_1$ to $z_1$.
   Obviously, $G_{20} \in{{\mathcal{U}}}_{2}(n,k,3)$, 
   where the cycle contained in $G_{20}$ is denoted by
   $C'_3 = u_1 \widetilde{e}_1 u_2 \widetilde{e}_2 z_1 g'' u_1$, each vertex of $\widetilde{e}_1 \setminus\{u_{1},u_{2}\}$ of $C'_3$ of $G_{20}$ is a core vertex
   and $\widetilde{e}_1$ is  the  perfect matching edge on $C'_3$ of $G_{20}$.
   By Lemma \ref{lemma2.1},  we obtain $\rho(G_{20})>\rho(G_{19})$.
   By Corollary \ref{corollary3.2}, we obtain $\max\{\rho(A_{n,k}),\rho(D_{n,k})\}\geq\rho(G_{20})$
   for $3\leq m\leq8$ and $\rho(D_{n,k})\geq\rho(G_{20})$ for $m\geq 9$, where $\rho(D_{n,k})=\rho(G_{20})$
   iff $G_{20}\cong D_{n,k}$.
   Therefore, we get $\rho(L_{n,k})>\rho(G_{20})$
   since $\rho(L_{n,k})>\rho(A_{n,k})$ and $\rho(L_{n,k})>\rho(D_{n,k})$ (by Lemma \ref{lemma444}).
   By the above proofs, we obtain $\rho(L_{n,k})>\rho(G)$ for $G\in\bar{\Gamma}_{1,1}(n,k)$.  Therefore, we get Lemma \ref{lemma44.5} in Case (i).

 \noindent\textbf{Case (ii).} $x_{u_{1}}\geq x_{z_{1}}$.

  Let $G_{21}$ be the hypergraph obtained from $G_{19}$ by moving $g'$ from $z_1$ to $u_1$ in such a way that
  $G_{21}$ has a perfect matching.
  Obviously, $G_{21}\in \bar{\Gamma}_{1,1}(n,k)$.
  By Lemma \ref{lemma2.1}, we obtain $\rho(G_{21})>\rho(G_{19})$.
  If $G_{21}\cong I_{n,k}$,
  then by Lemma \ref{lemma44.4} and  the above proofs, we get
  $\rho(L_{n,k})>\rho(I_{n,k})> \rho(G)$
   for $G\in\bar{\Gamma}_{1,1}(n,k)$. Namely, Lemma \ref{lemma44.5} holds.
  Otherwise, we assume $G_{21}\ncong I_{n,k}$.
  By repeatedly using the same procedure as those in Cases (i) and (ii),
  we finally obtain $\max\{\rho(L_{n,k}),\rho(I_{n,k})\}>\rho(G_{21})$.
  Therefore,  it follows from  $\rho(L_{n,k})>\rho(I_{n,k})$ (by Lemma \ref{lemma44.4}) and the above proofs, we have  Lemma \ref{lemma44.5}.
   ~~$\Box$

  By the methods similar to those for the proofs of Lemma \ref{lemma33.1}, we get
  Lemma \ref{lemma44.2}
  as follows.

 \begin{lemma}\label{lemma44.2}  
 Let $G\in\bar{\Gamma}_{1,2}(n,k)$, where  $n\geq k(k-1)$ with $k\geq3$.
  We have  $\rho(J_{n,k})\geq\rho(G)$ with the equality iff $G\cong J_{n,k}$.
 \end{lemma}

  \begin{lemma}\label{lemma44.3} 
 We have $\rho(I_{n,k})> \rho(J_{n,k})$ for  $n\geq 2k(k-1)$ and $k \geq 3$, where $I_{n,k}$ and $J_{n,k}$ are shown in Fig.\ \ref{fig-a}(b) and Fig.\ \ref{fig-5}(a), respectively.
   \end{lemma}

 \noindent\textbf{Proof}: \textcolor{black}{When $m\geq 2$, obviously} $I_{n,k}\ncong J_{n,k}$.

   We construct a weighted incidence matrix $B_{I_{n,k}}$ for $I_{n,k}$ as follows.
 Let $B_{I_{n,k}}(v,e)=0$ if $v\not\in e$;
   $B_{I_{n,k}}(v,e)=1$ if $v\in e$ and $v$ is a core vertex;
  $B_{I_{n,k}}(v,e)=\alpha$ if $v\in e$, $e$ is a pendent edge and $d_{I_{n,k}}(v)=2$;
 $B_{I_{n,k}}(v,e)=1-\alpha$ if $v\in e$ ($v\neq u_2$), $e$ is not a pendent edge and $d_{I_{n,k}}(v)=2$;
 $B_{I_{n,k}}(u_1,g_i)=\dfrac{\alpha}{(1-\alpha)^{k-1}}$ for $i=1,\cdots ,m-1$ and $m\geq 2$;
 $B_{I_{n,k}}(u_1,\widetilde{e}_{2})=x_{6}$;
 $B_{I_{n,k}}(u_2,\widetilde{e}_{2})=y_{6}$;
 $B_{I_{n,k}}(u_1,\widetilde{e}_{1})=c_{3}$; and
 $B_{I_{n,k}}(u_2,\widetilde{e}_{1})=d_{3}$,
%
  where $x_{6}$, $y_{6}$, $c_{3}$,  $d_{3}$, and $\alpha$ satisfy  (\ref{1.6})--(\ref{1.9}) as follows:
 \begin{numcases}{}
 x_{6}+c_{3}+\dfrac{(m-1)\alpha  }{(1-\alpha  )^{k-1}}=1,&\label{1.6}\\
 y_{6}+d_{3}=1,&\label{1.7}\\
 c_{3}d_{3}=\alpha  ,&\label{1.8}\\
 x_{6}y_{6}(1-\alpha  )^{k-2}=\alpha  ,&\label{1.9}\\
 \dfrac{x_{6}d_{3}}{y_{6}c_{3}}=1.\nonumber
 \end{numcases}
 We can check that  $\sum_{e:e\in E_{I_{n,k}}(v)}B_{I_{n,k}}(v,e)=1$ for any $v\in V(I_{n,k})$,
  $\prod_{v:v\in e}B_{I_{n,k}}(v,e)=\alpha$ for any $e\in E(I_{n,k})$, and $B_{I_{n,k}}$ is consistent.
   Thus, $I_{n,k}$ is consistently $\alpha$-normal.
    By Lemma \ref{lemma4}(i), we have $\rho({I_{n,k}})=\alpha^{-\frac{1}{k}}$.

  We construct a weighted incidence matrix $B_{J_{n,k}}$ for $J_{n,k}$ as follows.
 Let $B_{J_{n,k}}(v,e)=0$ if $v\not\in e$;
   $B_{J_{n,k}}(v,e)=1$ if $v\in e$ and $v$ is a core vertex;
  $B_{J_{n,k}}(v,e)=\alpha$ if $v\in e$, $e$ is a pendent edge and $d_{J_{n,k}}(v)=2$;
 $B_{J_{n,k}}(v,e)=1-\alpha$ if $v\in e$ ($v\neq u_1,u_2$), $e$ is not a pendent edge and $d_{J_{n,k}}(v)=2$;
 $B_{J_{n,k}}(u_{2,1},g_i)=\dfrac{\alpha}{(1-\alpha)^{k-1}}$ for $i=1,\cdots,m-1$ and $m\geq 2$;
 $B_{J_{n,k}}(u_2,\widetilde{e}_{2})=x_{7}$;
 $B_{J_{n,k}}(u_1,\widetilde{e}_{2})=y_{7}$;
 $B_{J_{n,k}}(u_2,\widetilde{e}_{1})=c_{3}$;
 $B_{J_{n,k}}(u_1,\widetilde{e}_{1})=d_{3}$; and
 $B_{J_{n,k}}(u_{2,1},\widetilde{e}_{2})=1-\alpha  -\dfrac{(m-1)\alpha  }{(1-\alpha  )^{k-1}}$,
 %
 where $x_{7}$, $y_{7}$, $c_{3}$, and $d_{3}$ satisfy (\ref{1.10}) and (\ref{1.11}) as follows:
 \begin{numcases}{}
 x_{7}+c_{3}=1,&\label{1.10} \\
 y_{7}+d_{3}=1.&\label{1.11}
 \end{numcases}
  We can check that  $\sum_{e:e\in E_{J_{n,k}}(v)}B_{J_{n,k}}(v,e)=1$ for any $v\in V(J_{n,k})$ and
  $\prod_{v:v\in e}B_{J_{n,k}}(v,e)=\alpha$ for  any $e\in E(J_{n,k})$ and $e\neq \widetilde{e}_{2}$.
  Next, we prove $\prod_{v: v\in \widetilde{e}_{2}}B_{J_{n,k}}(v,\widetilde{e}_{2})>\alpha$.

   Let $\dfrac{x_{6}}{y_{6}}=\dfrac{c_{3}}{d_{3}}=b$.
   We have $x_{6}=by_{6}$ and $c_{3}=bd_{3}$.
      Substituting $c_{3}=bd_{3}$  into (\ref{1.8}), we get $\alpha=bd^2_{3}$.
  Substituting $x_{6}=by_{6}$, $\alpha=bd^2_{3}$ and (\ref{1.7}) into (\ref{1.9}), we get
   $d_{3}=\dfrac{(1-\alpha  )^{k/2-1}}{1+(1-\alpha  )^{k/2-1}}$.
   Therefore, by (\ref{1.7}), we have
   \begin{align}\label{66}
    y_{6}=1-d_{3}=\dfrac{1}{1+(1-\alpha  )^{k/2-1}}.
    \end{align}
    Substituting $x_{6}=by_{6}$
    and $c_{3}=bd_{3}$ into (\ref{1.6}) and bearing (\ref{1.7}) (namely $y_{6}+d_{3}=1$),  $\alpha=bd^2_{3}$ and  $d_{3}=\dfrac{(1-\alpha  )^{k/2-1}}{1+(1-\alpha  )^{k/2-1}}$ in mind,
   we get
    \begin{align}\label{55mm}
    \dfrac{(m-1)\alpha  }{\left(1-\alpha  \right)^{k-1}}=1-b\left(y_{6}+d_{3}\right)=1-b=1-\alpha \dfrac{\left(1+\left(1-\alpha  \right)^{k/2-1}\right)^2}{\left(1-\alpha  \right)^{k-2}}.
    \end{align}
 We obtain
 \begin{align}
  &\prod_{v: v\in \widetilde{e}_{2}}B_{J_{n,k}}(v,\widetilde{e}_{2})-\alpha  \nonumber\\
  &=x_{7}y_{7}\left[1-\alpha  -\dfrac{(m-1)\alpha  }{(1-\alpha  )^{k-1}}\right](1-\alpha  )^{k-3}-\alpha  \nonumber\\
   &=y_{6}\left[x_{6}+\dfrac{(m-1)\alpha  }{(1-\alpha  )^{k-1}}\right]\left[1-\alpha  -\dfrac{(m-1)\alpha  }{(1-\alpha  )^{k-1}}\right](1-\alpha  )^{k-3}-\alpha  \label{11}\\
  &=\left[\dfrac{\alpha  }{(1-\alpha  )^{k-2}}+\dfrac{(m-1)\alpha  }{(1-\alpha  )^{k-1}}y_{6}\right]\left[1-\alpha  -\dfrac{(m-1)\alpha  }{(1-\alpha  )^{k-1}}\right](1-\alpha  )^{k-3}-\alpha  \label{22}\\
  &=\dfrac{(m-1)\alpha  }{1-\alpha  }\left[y_{6}\left(1-\dfrac{(m-1)\alpha  }{(1-\alpha  )^{k-1}}\cdot\frac{1}{1-\alpha  }\right)-\dfrac{\alpha  }{(1-\alpha  )^{k-1}}\right] \label{w22}\\
  &=\dfrac{(m-1)\alpha  ^2}{(1-\alpha  )^2}\left[\dfrac{1}{(1-\alpha  )^{k/2-1}}-\dfrac{1}{1+(1-\alpha  )^{k/2-1}}\right] \label{2ss2}\\
  &>0.\label{77}
     \end{align}
 It is noted that (\ref{11}) follows from $x_{7}=x_{6}+\dfrac{(m-1)\alpha  }{(1-\alpha  )^{k-1}}$
 (by (\ref{1.6}) and (\ref{1.10})) and $y_{7}=y_{6}$ (by (\ref{1.7}) and (\ref{1.11})).
 Substituting $x_{6}y_{6}=\alpha/(1-\alpha  )^{k-2}$ (by (\ref{1.9})) into (\ref{11}), we have (\ref{22}).
 By calculation, we get  (\ref{w22}).
 Substituting the expression of $y_{6}$ (namely (\ref{66})) and the expression of $\dfrac{(m-1)\alpha  }{(1-\alpha  )^{k-1}}$  (namely (\ref{55mm})) into (\ref{w22}), we have (\ref{2ss2}).
 Since $m>1$ and $0<\alpha  <1$, we get (\ref{77}).

 By the above proofs, we get that   $J_{n,k}$ is strictly $\alpha$-subnormal.
 Therefore, by Lemma \ref{lemma4}(ii), we have $\rho({J_{n,k}})<\alpha^{-\frac{1}{k}}$.
 Thus, by Lemma \ref{lemma4},  we obtain $\rho(I_{n,k})>\rho(J_{n,k})$ for  $n\geq 2k(k-1)$ and $k \geq 3$.
 ~~$\Box$

 The  hypergraph with the maximal spectral radius in $\Gamma_{1}(n,k)$ are shown in Corollary \ref{corollary4.2}.

  \begin{corollary}\label{corollary4.2}
 Let $G\in\Gamma_{1}(n,k)$, where $n\geq 2k(k-1)$ with $k \geq 3$.
  We have $\rho(L_{n,k})>\rho(G)$.
 \end{corollary}

 \noindent\textbf{Proof}: Let $n\geq 2k(k-1)$ with $k\geq3$.
 By Lemma \ref{lemma44.5}, we have
 $\rho(L_{n,k})>\rho(G)$ for  $G\in\bar{\Gamma}_{1,1}(n,k)$.
 By Lemmas \ref{lemma44.4}, \ref{lemma44.2} and \ref{lemma44.3}, we obtain
 $\rho(L_{n,k})>\rho(I_{n,k})>\rho(J_{n,k})>\rho(G)$ for  $G\in\bar{\Gamma}_{1,2}(n,k)$.
 Therefore, we have  $\rho(L_{n,k})>\rho(G)$
 for  $G\in \bar{\Gamma}_{1}(n,k)$ since $\bar{\Gamma}_{1,1}(n,k)\cup \bar{\Gamma}_{1,2}(n,k)=\bar{\Gamma}_{1}(n,k)$.
  Furthermore, by Lemma \ref{lemma4.3}, we get Corollary \ref{corollary4.2}.
  ~~$\Box$


%


 \begin{lemma}\label{lemma4.4}
  Let $G\in{\Gamma}_{2}(n,k)$,  where $n\geq 2k(k-1)$ and $k \geq3$.
   There exists one hypergraph $\overline{G}\in\bar{\Gamma}_{2}(n,k)$ such that $\rho(\overline{G})\geq\rho(G)$ with the equality iff $G\cong \overline{G}$.
 \end{lemma}

 \noindent\textbf{Proof}: Let $G\in {\Gamma}_{2}(n,k)$ with $n\geq 2k(k-1)$ and $k \geq3$.
  By 
  Lemma \ref{lemmawang1},
  we get that there exists one hypergraph $\dot{G}$ such that $\rho(\dot{G})\geq\rho(G)$ with the equality iff $G\cong \dot{G}$,
  where $\dot{G}$ satisfies two conditions: (i) each vertex in $C_2$ must be attached by a pendent edge, and
 (ii) the vertex in  $C_2$  maybe  attached by a hypertree which has at least $k$ edges, where $k\geq 3$.

  If $\dot{G}\in\bar{\Gamma}_{2}(n,k)$, then  Lemma \ref{lemma4.4} holds.
   Next, suppose  $\dot{G}\notin\bar{\Gamma}_{2}(n,k)$.
   In $\dot{G}$, there exist at least one vertex in $\widetilde{e}_{1}\cup \widetilde{e}_{2}$ which is attached by a hypertree having at least $k\geq 3$ edges.
  In $\dot{G}$, let $V_1(\dot{G})$ be the subset of $V(\dot{G})$ in which each vertex is
  attached by a hypertree having at least $k\geq 3$ edges, where $|V_1(\dot{G})| \geq 1$.
    Let $w$ be  a vertex in $V_1(\dot{G})$ and let $w\in \widetilde{e}_{1}$.
   Let  $\bm{x}$ be the principal eigenvector of $\dot{G}$ corresponding to $\rho(\dot{G})$.
   Among all the vertices in $V_1(\dot{G})$,  we suppose that $w$ has the maximal component $x_w$  among  $\bm{x}$.
By the methods similar to those for the second paragraph, Cases (i) and (ii) in Lemma \ref{lemma3.4}, we obtain a hypergraph $\overline{G}\in\bar{\Gamma}_{2}(n,k)$ such that
$\rho(\overline{G})>\rho(\dot{G}) \geq\rho(G)$. Therefore, we have  Lemma \ref{lemma4.4}.
  ~~$\Box$
%
%

%

%




  By the methods similar to those for the proofs of Lemma \ref{lemma33.1}, we get Lemma \ref{lemma44.1}   as follows.

 \begin{lemma}\label{lemma44.1}
 Let $G\in\bar{\Gamma}_{2}(n,k)$, where  $n\geq 2k(k-1)$ with $k\geq3$.
 We have  $\rho(L_{n,k})\geq\rho(G)$ with the equality iff $G\cong L_{n,k}$.
 \end{lemma}

  By   Lemmas \ref{lemma4.4} and \ref{lemma44.1},
 we have Corollary \ref{corollary4.1}  as follows.

 \begin{corollary}\label{corollary4.1}
 Let $G\in\Gamma_{2}(n,k)$, where $n\geq 2k(k-1)$ with $k\geq 3$.
 We have $\rho(L_{n,k})\geq\rho(G)$ with the equality iff $G\cong L_{n,k}$.
 \end{corollary}

  By Corollaries \ref{corollary4.2} and \ref{corollary4.1}, 
    we get Theorem \ref{theorem4.6}.

 \begin{theorem}\label{theorem4.6}
 Let $G\in\Gamma(n,k)$, where $n=mk(k-1)$ with $m\geq1$ and $k\geq 3$.
 When $m=1$, $G\cong I_{n,k}\cong J_{n,k}$ and $\rho(G)=\rho(I_{n,k})=\rho(J_{n,k})$.
 When $m\geq 2$, $\rho(L_{n,k})\geq\rho(G)$ with the equality iff $G\cong L_{n,k}$.
 \end{theorem}

 \section{The hypergraph with the maximal spectral radius among  $\mathcal{U}(n,k)\cup\Gamma(n,k)$}

  In  Section 5, we get the hypergraph with the maximal spectral radius among  $\mathcal{U}(n,k)\cup\Gamma(n,k)$, where $n=mk(k-1)$, $m\geq 1$ and $k\geq 3$. 
   For $n=k(k-1)$ with $k\geq 3$,  there is only
  one hypergraph  $I_{n,k}$  (namely $J_{n,k}$) among  $\mathcal{U}(n,k)\cup\Gamma(n,k)$.
  For $n=mk(k-1)$  with $m\geq 2$ and $k\geq 3$,  
   we have Theorem \ref{theorem5.2} as follows.

  \begin{theorem}\label{theorem5.2}
 Let $G\in \mathcal{U}(n,k)\cup\Gamma(n,k)$ with  $n\geq2k(k-1)$ and $k\geq 3$,
  we have $\rho(L_{n,k})\geq\rho(G)$ with the equality iff $G\cong L_{n,k}$.
 \end{theorem}

 \noindent\textbf{Proof}: Let $n\geq2k(k-1)$ and $k\geq 3$.
   By Theorem \ref{theorem3.6} and Lemma  \ref{lemma444}, we have $\rho(L_{n,k})>\rho(G)$ for
 $G\in \mathcal{U}(n,k)$.
  By Theorem \ref{theorem4.6}, we have $\rho(L_{n,k})\geq\rho(G)$ for
 $G\in \Gamma(n,k)$ with the equality iff $G\cong L_{n,k}$.
 Therefore, we get  Theorem \ref{theorem5.2}. ~~$\Box$

\bibliographystyle{elsevier_citation_order}  
\bibliography{Perfect_matching_unicyclic_hypergraph}


\newpage

\begin{figure}
\thicklines
\begin{tikzpicture}
\draw[thick,black](0.8,0) ellipse(1.5 and 0.2); 
\draw[thick,black](3,0)  ellipse(1.5 and 0.2); 
\fill(-0.35,0) circle(0.02);  
\fill(0.6,0) circle(0.02);  
\fill(0.40,0) circle(0.02);  
\fill(0.2,0) circle(0.02);  
\fill(1.2,0) circle(0.02); 
\draw[black](1.7,0)  circle(0.02);  
\draw[black](2.05,0)  circle(0.02);  
\fill(2.7,0) circle(0.02);   
\fill(3.5,0) circle(0.02);   
\fill(3.3,0) circle(0.02);   
\fill(3.1,0) circle(0.02);   
\fill(4,0) circle(0.02);   
\draw(-1,0) node{{\scriptsize $\widetilde{e}_{1}$}};
\draw(4.75,0) node{{\scriptsize $\widetilde{e}_{2}$}};
\draw(2.2,-0.4) node{{\scriptsize $u_2$}};
\draw(1.6,-0.4) node{{\scriptsize $u_{1}$}};
\draw(1.2,-0.55) node{{\scriptsize $u_{1,1}$}};
\draw(-0.35,-0.5) node{{\scriptsize $u_{1,k-2}$}};
\draw(2.70,-0.55) node{{\scriptsize $u_{2,1}$}};
\draw(4,-0.5) node{{\scriptsize $u_{2,k-2}$}};
\draw(1.9,-2) node{{\scriptsize $(a)$}};
\end{tikzpicture}
\begin{tikzpicture}
\draw[thick,black,rotate=45](2.4,-1.2) ellipse(1.3 and 0.2); 
\draw[thick,black,rotate=135](-0.05,-1.2) ellipse(1.3 and 0.2); 
\draw[thick,blue](-0.42,1.2) ellipse(1.0 and 0.15); 
\draw[black](0.42,1.2) circle(0.02);  
\draw[thick,blue](3.7,1.2)  ellipse(1.0 and 0.15); 
\draw[black](2.9,1.2) circle(0.02);  
\draw[thick,blue](0.19,0.55) ellipse(1.0 and 0.15); 
\draw[black](1.1,0.55) circle(0.02);  
\draw[thick,blue](3.1,0.55)  ellipse(1.0 and 0.15); 
\draw[black](2.22,0.55) circle(0.02);  
\draw[thick,black](0.8,0) ellipse(1.5 and 0.2); 
\draw[thick,red](3,0)  ellipse(1.5 and 0.2); 
\fill(1.0,1.8) circle(0.5pt);
\fill(1.5,1.8) circle(0.5pt);
\fill(2.0, 1.8) circle(0.5pt); 
\fill(0.8, 1) circle(0.5pt);
\fill(0.88, 0.9) circle(0.5pt);
\fill(0.98, 0.8) circle(0.5pt);
\fill(2.4, 0.8) circle(0.5pt);
\fill(2.5, 0.9) circle(0.5pt);
\fill(2.6, 1.0) circle(0.5pt);
\draw[thick,blue,rotate=270](0.9,-0.35) ellipse(1 and 0.15); 
\draw[thick,blue,rotate=270](0.9,1) ellipse(1 and 0.15); 
\draw[black](-0.35,0) circle(0.02);  
\fill(0.40,0) circle(0.02);  
\fill(0.01,0) circle(0.02);  
\fill(0.75,0) circle(0.02); 
\fill(0.40,-0.7) circle(0.02);  
\fill(0.01,-0.7) circle(0.02);  
\fill(0.75,-0.7) circle(0.02); 
\draw[black](1.0,0)  circle(0.02);  
\draw[black](1.7,0)  circle(0.02);  
\draw[black](2.05,0)  circle(0.02);  
\fill(2.7,0) circle(0.02);   
\fill(3.3,0) circle(0.02);   
\fill(4,0) circle(0.02);   
\draw(1.5,2.5) node{{\scriptsize $m-1$}};
 \draw(1.63,2) node{$\overbrace{~~~~~~~~~~~~~~~~~~~}$};
\draw(-0.2,2) node{{\scriptsize $g_1$}};
\draw(3.7,2) node{{\scriptsize $g_{m-1}$}};
\draw(-1,0) node{{\scriptsize $\widetilde{e}_{2}$}};
\draw(4.75,0) node{{\scriptsize $\widetilde{e}_{1}$}};
\draw(2.2,-0.4) node{{\scriptsize $u_2$}};
\draw(1.57,-0.4) node{{\scriptsize $u_1$}};
\draw(2,-3) node{{\scriptsize $(a)$}};
\end{tikzpicture}
\caption{\label{fig-a} (a) $C_{2}$ and  (b) $I_{n,k}$}
\end{figure}
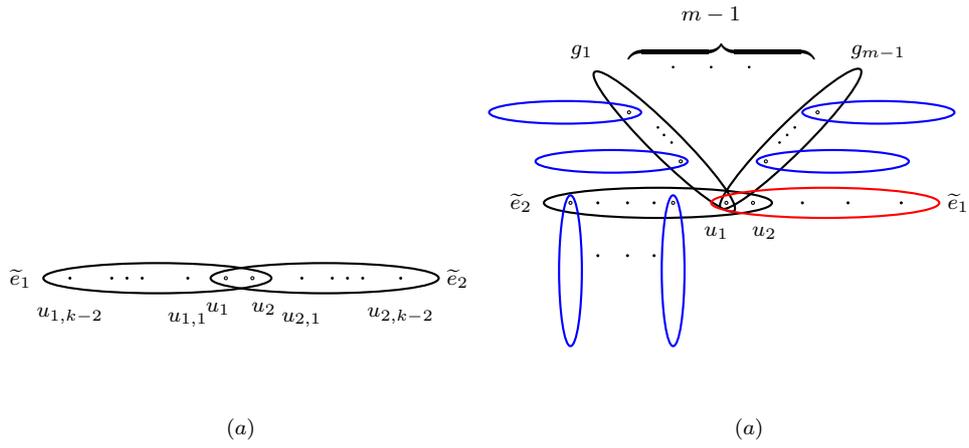

%
\end{document}